\documentclass[12pt]{elsarticle}
\usepackage[cp1250]{inputenc}
\usepackage{amsfonts}
\usepackage{amsmath}
\usepackage{graphicx}
\usepackage{epsfig}
\usepackage{amssymb}
\usepackage{amstext}
\usepackage{amsthm}
\usepackage{tkz-berge}

\usetikzlibrary{positioning}

\newtheorem{theorem}{Theorem}[section]
\newtheorem{lemma}{Lemma}[section]
\theoremstyle{definition}
\newtheorem{remark}{Remark}[section]
\newtheorem{corollary}{Corollary}[section]
\newtheorem{definition}{Definition}

\newenvironment{PrfFact}{{\bf Proof }}{{\hfill\tiny{$\blacksquare$\\}}}

\def\kob{[fill=white]circle (3pt)} 
\def\koc{[fill=black]circle (3pt)} 
\def\kwb{ [fill=white,minimum size=3pt] ++(-3pt,-3pt) rectangle ++(6pt,6pt)}
\def\kwc{ [fill=black,minimum size=3pt] ++(-3pt,-3pt) rectangle ++(6pt,6pt)}

\def\nowyii{
\newbox\AA
\savebox\AA{
\begin{tikzpicture}[
]

\tikzset{VertexStyle/.style = {
shape = circle,
fill = black,
inner sep = 0pt,
outer sep = 0pt,
minimum size = 0pt,
draw}}

  \SetVertexNoLabel
\grPath[Math,prefix=p,RA=1,RS=0]{4}
\begin{scope}[xshift=-0.5 cm]
\grPath[Math,prefix=q,RA=1,RS=1]{5}
\end{scope}
\begin{scope}[xshift=1.5 cm]
\grPath[Math,prefix=r,RA=1,RS=2]{1}
\end{scope}
\begin{scope}[xshift=1.5 cm]
\grPath[Math,prefix=n,RA=1,RS=-1]{1}
\Edges(p0,q0,p0,q1,p1,q2,p2,q3,p3,q4)
\Edges(n0,p0,n0,p1,n0,p2,n0,p3)
\Edges(r0,q0,r0,q1,r0,q2,r0,q3,r0,q4)

\foreach \k in {p1,q0,q4} {  \draw (\k) \kwc  ;}
\foreach \k in {p2,r0} {  \draw (\k) \kwb  ;}
\foreach \k in {p0,p3,q2} {  \draw (\k) \koc  ;}
\foreach \k in {n0,q1,q3} {  \draw (\k) \kob  ;}

\draw (r0) node[above,yshift=3pt]{$b$};
\draw (q0) node[above,yshift=1pt]{$c_{p-3}$};
\draw (q2) node[above right]{$c_{p}$};
\draw (q4) node[above,yshift=1pt]{$d_2$};
\draw (p0) node[below,xshift=-3pt]{$c_{p-2}$};
\draw (p1) node[below right,xshift=1pt]{$c_{p-1}$};
\draw (p2) node[below right,xshift=1pt]{$c$};
\draw (p3) node[below right,xshift=-1pt]{$d_1$};
\draw (n0) node[below,yshift=-3pt]{$a$};

\end{scope}

\end{tikzpicture}
}

\newbox\AB
\savebox\AB{
\begin{tikzpicture}
\tikzset{VertexStyle/.style = {
shape = circle,
fill = black,
inner sep = 0pt,
outer sep = 0pt,
minimum size = 0pt,
draw}}

  \SetVertexNoLabel
 
\grPath[Math,prefix=p,RA=1,RS=0]{4}
\begin{scope}[xshift=-0.5 cm]
\grPath[Math,prefix=q,RA=1,RS=1]{5}
\end{scope}
\begin{scope}[xshift=1.5 cm]
\grPath[Math,prefix=r,RA=1,RS=2]{1}
\end{scope}
\begin{scope}[xshift=1.5 cm]
\grPath[Math,prefix=n,RA=1,RS=-1]{1}
\Edges(p0,q0,p0,q1,p1,q2,p2,q3,p3,q4)
\Edges(n0,p0,n0,p1,n0,p2,n0,p3)
\Edges(r0,q0,r0,q1,r0,q2,r0,q3,r0,q4)

\foreach \k in {p2,q0,q4} {  \draw (\k) \kwc  ;}
\foreach \k in {p1,r0} {  \draw (\k) \kwb  ;}
\foreach \k in {p0,p3,q2} {  \draw (\k) \koc  ;}
\foreach \k in {n0,q1,q3} {  \draw (\k) \kob  ;}

\draw (r0) node[above,yshift=3pt]{$b$};
\draw (q0) node[above,yshift=1pt]{$c_{p-3}$};
\draw (q2) node[above right]{$c_{p}$};
\draw (q4) node[above,yshift=1pt]{$d_2$};
\draw (p0) node[below,xshift=-3pt]{$c_{p-2}$};
\draw (p1) node[below right,xshift=1pt]{$c_{p-1}$};
\draw (p2) node[below right,xshift=1pt]{$c$};
\draw (p3) node[below right,xshift=-1pt]{$d_1$};
\draw (n0) node[below,yshift=-3pt]{$a$};

\end{scope}

\end{tikzpicture}
}

\begin{figure}[!htb]
\centering
\begin{tikzpicture}

\node (AA) {\usebox\AA };

\node[right=of AA] (AB) {
\usebox\AB}	;
\draw[->] (AA) -- (AB) node[pos=0.5,above]{$A'$};
\end{tikzpicture}
\caption{An edge $c_{p-1}c$ is a Kempe chain of $A(2, 4)$ and the  colouring $A^{'}$}
\end{figure}

} 

\def\nowyi{

\begin{figure}[!htb]
\centering
\begin{tikzpicture}[
koc/.style={circle,draw=black,fill=black,,inner sep=0pt, outer sep=0pt,thick,minimum size=6pt},
kob/.style={circle,draw=black,fill=white,thick,inner sep=0pt, outer sep=0pt,minimum size=6pt},
kwc/.style={rectangle,draw=black,fill=black,thick,inner sep=0pt, outer sep=0pt,minimum size=6pt},
kwb/.style={rectangle,draw=black,fill=white,thick,inner sep=0pt, outer sep=0pt,minimum size=6pt},
]

\node[kob] (v0) at (0:0) {};
\node[kwb] (a1) at (-360/14:1.2) {};
\node[koc] (a2) at (360/14:1.2) {};
\node[kwb] (a3) at (3*360/14:1.2) {};
\node[koc] (a4) at (5*360/14:1.2) {};
\node[kwb] (a5) at (7*360/14:1.2) {};
\node[koc] (a6) at (9*360/14:1.2) {};
\node[kwc] (a7) at (11*360/14:1.2) {};

\node[kwc] (b1) at (0:2.1) {};
\node[kob] (b2) at (360/7:2.1) {};
\node[kwc] (b3) at (2*360/7:2.1) {};
\node[kob] (b4) at (3*360/7:2.1) {};
\node[kwc] (b5) at (4*360/7:2.1) {};
\node[kob] (b6) at (5*360/7:2.1) {};
\node[koc] (b7) at (6*360/7:2.1) {};

\foreach \k in {1,2,3,4,5,6,7} {\draw (v0) -- (a\k);
\draw (a\k) -- (b\k);
\draw (b3) -- (a4) node[pos=0.7,  above, xshift=-3pt]{$e$};
}

\foreach \k in {1,2,3,4,5,6}
{\draw (a\k) -- (a\the\numexpr\k+1\relax);
\draw (b\k) -- (a\the\numexpr\k+1\relax);
\draw (b\k) -- (b\the\numexpr\k+1\relax);
}
\draw (a7) -- (a1);
\draw (b7) -- (a1);
\draw (b7) -- (b1);
\end{tikzpicture}

\caption{The graph $H_{7} = G_{7} - b$ and the colouring $Q_{2,e}$ restricted to $H_{7}$}
\end{figure}
} 

\def\nowyiii{
\newbox\piel
\savebox\piel{
\begin{tikzpicture}
\tikzset{VertexStyle/.style = {
shape = circle,
fill = black,
inner sep = 0pt,
outer sep = 0pt,
minimum size = 0pt,
draw}}

\SetVertexNoLabel

\grPath[Math,prefix=p,RA=1,RS=0]{4}
\begin{scope}[xshift=-0.5 cm]
\grPath[Math,prefix=q,RA=1,RS=1]{5}
\end{scope}

\Edges(q0,p0,q1,p1,q2,p2,q3,p3,q4)
\Edges(q0,q4)
\Edges(p0,p3)

\foreach \k in {q2} {  \draw (\k) \kwb  ;}
\foreach \k in {p0,p2,q4} {  \draw (\k) \kwc  ;}
\foreach \k in {q1,q3} {  \draw (\k) \kob  ;}
\foreach \k in {p1,p3,q0} {  \draw (\k) \koc  ;}

\draw (p2) node[below,yshift=-1pt]{$y$};
\draw (q1) node[above,yshift=1pt]{$x_1$};
\draw (q2) node[above,yshift=1pt]{$x_2$};
\draw (q3) node[above,yshift=1pt]{$x_3$};

\end{tikzpicture}
}

\newbox\pier
\savebox\pier{
\begin{tikzpicture}
\tikzset{VertexStyle/.style = {
shape = circle,
fill = black,
inner sep = 0pt,
outer sep = 0pt,
minimum size = 0pt,
draw}}

\SetVertexNoLabel

\grPath[Math,prefix=p,RA=1,RS=0]{4}
\begin{scope}[xshift=-0.5 cm]
\grPath[Math,prefix=q,RA=1,RS=1]{5}
\end{scope}

\Edges(q0,p0,q1,p1,q2,p2,q3,p3,q4)
\Edges(q0,q4)
\Edges(p0,p3)

\foreach \k in {p2} {  \draw (\k) \kwb  ;}
\foreach \k in {p0,q2,q4} {  \draw (\k) \kwc  ;}
\foreach \k in {q1,q3} {  \draw (\k) \kob  ;}
\foreach \k in {p1,p3,q0} {  \draw (\k) \koc  ;}

\draw (p2) node[below,yshift=-1pt]{$\phantom{y}$};
\draw (q1) node[above,yshift=1pt]{$x_1$};
\draw (q2) node[above,yshift=1pt]{$x_2$};
\draw (q3) node[above,yshift=1pt]{$x_3$};

\end{tikzpicture}
}

\begin{figure}[!htb]
\centering
\begin{tikzpicture}
\node (piel) {\usebox\piel};
\node[right=of piel] (pier) {\usebox\pier};

\draw[->] (piel) -- (pier) node[pos=.5,above] {$B_2$};

\end{tikzpicture}
\caption{A Kempe chain $x_{2}y$ of $B_{1}(2, 4)$ and the colouring $B_{2}$.}
\end{figure}
} 

\def\nowyiv{
\newbox\drul
\savebox\drul{
\begin{tikzpicture}
\tikzset{VertexStyle/.style = {
shape = circle,
fill = black,
inner sep = 0pt,
outer sep = 0pt,
minimum size = 0pt,
draw}}

\SetVertexNoLabel

\grPath[Math,prefix=p,RA=1,RS=0]{6}
\begin{scope}[xshift=0.5 cm]
\grPath[Math,prefix=q,RA=1,RS=1]{5}
\end{scope}

\Edges(p0,q0,p1,q1,p2,q2,p3,q3,p4,q4,p5)
\Edges(q0,q4)
\Edges(p0,p5)

\foreach \k in {p1,p4,q2} {  \draw (\k) \kwb  ;}
\foreach \k in {p3,p5,q0} {  \draw (\k) \kwc  ;}
\foreach \k in {q1,q3} {  \draw (\k) \kob  ;}
\foreach \k in {p0,p2,q4} {  \draw (\k) \koc  ;}

\draw (p1) node[below,yshift=-1pt]{$x_1$};
\draw (p2) node[below,yshift=-1pt]{$y$};
\draw (p3) node[below,yshift=-1pt]{$z$};
\draw (p4) node[below,yshift=-1pt]{$x_5$};
\draw (q1) node[above,yshift=1pt]{$x_2$};
\draw (q2) node[above,yshift=1pt]{$x_3$};
\draw (q3) node[above,yshift=1pt]{$x_4$};

\end{tikzpicture}
}

\newbox\drur
\savebox\drur{
\begin{tikzpicture}
\tikzset{VertexStyle/.style = {
shape = circle,
fill = black,
inner sep = 0pt,
outer sep = 0pt,
minimum size = 0pt,
draw}}

\SetVertexNoLabel

\grPath[Math,prefix=p,RA=1,RS=0]{6}
\begin{scope}[xshift=0.5 cm]
\grPath[Math,prefix=q,RA=1,RS=1]{5}
\end{scope}

\Edges(p0,q0,p1,q1,p2,q2,p3,q3,p4,q4,p5)
\Edges(q0,q4)
\Edges(p0,p5)

\foreach \k in {p1,p4,q2} {  \draw (\k) \kwb  ;}
\foreach \k in {p2,p5,q0} {  \draw (\k) \kwc  ;}
\foreach \k in {q1,q3} {  \draw (\k) \kob  ;}
\foreach \k in {p0,p3,q4} {  \draw (\k) \koc  ;}

\draw (p1) node[below,yshift=-1pt]{$\phantom{x_1}$};
\draw (q1) node[above,yshift=1pt]{$x_2$};
\draw (q3) node[above,yshift=1pt]{$x_4$};

\end{tikzpicture}
}

\begin{figure}[!htb]
\centering
\begin{tikzpicture}
\node (drul) {\usebox\drul};
\node[right=of drul] (drur) {\usebox\drur};

\draw[->] (drul) -- (drur) node[pos=.5,above] {$B_1$};

\end{tikzpicture}
\caption{A Kempe chain $yz$ of $A(3, 4)$ and the colouring  $B_1$}
\end{figure}
} 

\def\nowyv{
\newbox\trzl
\savebox\trzl{
\begin{tikzpicture}
\tikzset{VertexStyle/.style = {
shape = circle,
fill = black,
inner sep = 0pt,
outer sep = 0pt,
minimum size = 0pt,
draw}}

\SetVertexNoLabel

\grPath[Math,prefix=p,RA=1,RS=0]{4}
\begin{scope}[xshift=-0.5 cm]
\grPath[Math,prefix=q,RA=1,RS=1]{4}
\end{scope}

\Edges(q0,p0,q1,p1,q2,p2,q3,p3)
\Edges(q0,q3)
\Edges(p0,p3)

\foreach \k in {p3,q2} {  \draw (\k) \kwb  ;}
\foreach \k in {p1,q0} {  \draw (\k) \kwc  ;}
\foreach \k in {q1,q3} {  \draw (\k) \kob  ;}
\foreach \k in {p0,p2} {  \draw (\k) \koc  ;}

\draw (p1) node[below,yshift=-1pt]{$y$};
\draw (p3) node[below,yshift=-1pt]{$x_4$};
\draw (q1) node[above,yshift=1pt]{$x_1$};
\draw (q2) node[above,yshift=1pt]{$x_2$};
\draw (q3) node[above,yshift=1pt]{$x_3$};

\end{tikzpicture}
}

\newbox\trzr
\savebox\trzr{
\begin{tikzpicture}
\tikzset{VertexStyle/.style = {
shape = circle,
fill = black,
inner sep = 0pt,
outer sep = 0pt,
minimum size = 0pt,
draw}}

\SetVertexNoLabel

\grPath[Math,prefix=p,RA=1,RS=0]{4}
\begin{scope}[xshift=-0.5 cm]
\grPath[Math,prefix=q,RA=1,RS=1]{4}
\end{scope}

\Edges(q0,p0,q1,p1,q2,p2,q3,p3)
\Edges(q0,q3)
\Edges(p0,p3)

\foreach \k in {p1,p3} {  \draw (\k) \kwb  ;}
\foreach \k in {p0,q2} {  \draw (\k) \kwc  ;}
\foreach \k in {q1,q3} {  \draw (\k) \kob  ;}
\foreach \k in {p2,q0} {  \draw (\k) \koc  ;}

\draw (p3) node[below,yshift=-1pt]{$\phantom{x_4}$};
\draw (q1) node[above,yshift=1pt]{$\phantom{x_1}$};

\end{tikzpicture}
}

\begin{figure}[!htb]
\centering
\begin{tikzpicture}
\node (trzl) {\usebox\trzl};
\node[right=of trzl] (trzr) {\usebox\trzr};

\draw[->] (trzl) -- (trzr) node[pos=.5,above] {$B$};

\end{tikzpicture}
\caption{A Kempe chain $x_{2}y$  of $A(2, 4)$ and the colouring $B$}
\end{figure}
} 

\def\nowyvi{
\begin{figure}[!htb]
\centering

\begin{tikzpicture}
\tikzset{VertexStyle/.style = {
shape = circle,
fill = black,
inner sep = 0pt,
outer sep = 0pt,
minimum size = 0pt,
draw}}

\SetVertexNoLabel

\grEmptyPath[Math,prefix=p,RA=1,RS=0]{10}
\begin{scope}[xshift=-0.5 cm]
\grEmptyPath[Math,prefix=q,RA=1,RS=1]{11}
\end{scope}

\Edges(q0,p0,q1,p1,q2,p2,q3)
\Edges(q0,q3)
\Edges(p0,p2)
\Edges(q7,p7,q8,p8,q9,p9,q10)
\Edges(q7,q10)
\Edges(p7,p9)
\Edges[style={dotted}](q3,q7)
\Edges[style={dotted}](p2,p7)

\foreach \k in {p2,q2,q8} {  \draw (\k) \kwb  ;}
\foreach \k in {p0,p2,p8,q7,q10} {  \draw (\k) \kwc ;}
\foreach \k in {q1,q9} {  \draw (\k) \kob  ;}
\foreach \k in {p1,p7,p9,q0,q3} {  \draw (\k) \koc  ;}

\draw (p0) node[below,yshift=-1pt]{$a_2$};
\draw (p1) node[below,yshift=-1pt]{$a_3$};
\draw (p2) node[below,yshift=-1pt]{$a_4$};

\draw (p7) node[below,yshift=-1pt]{$b_2$};
\draw (p8) node[below,yshift=-1pt]{$b_3$};
\draw (p9) node[below,yshift=-1pt]{$b_4$};

\draw (q0) node[above,yshift=1pt]{${a_1}$};
\draw (q1) node[above,yshift=1pt]{${x_1}$};
\draw (q2) node[above,yshift=1pt]{${x_2}$};
\draw (q3) node[above,yshift=1pt]{$a_5$};

\draw (q7) node[above,yshift=1pt]{$b_1$};
\draw (q8) node[above,yshift=1pt]{$y_1$};
\draw (q9) node[above,yshift=1pt]{$y_2$};
\draw (q10) node[above,yshift=1pt]{$b_5$};

\draw (q5) node[above,yshift=1pt]{$N_{b}$};
\draw (p5) node[above,xshift=-5mm]{$N_{a}$};
\end{tikzpicture}
\caption{Vertices $x_{1}$ and $y_2$ are coloured the same by $B_{1}$}
\end{figure}
} 

\def\nowyvii{
\begin{figure}[!htb]
\centering

\begin{tikzpicture}
\tikzset{VertexStyle/.style = {
shape = circle,
fill = black,
inner sep = 0pt,
outer sep = 0pt,
minimum size = 0pt,
draw}}

\SetVertexNoLabel

\grEmptyPath[Math,prefix=p,RA=1,RS=0]{10}
\begin{scope}[xshift=-0.5 cm]
\grEmptyPath[Math,prefix=q,RA=1,RS=1]{11}
\end{scope}

\Edges(q0,p0,q1,p1,q2,p2,q3)
\Edges(q0,q3)
\Edges(p0,p2)
\Edges(p7,q8,p8,q9,p9,q10)
\Edges(q8,q10)
\Edges(p7,p9)
\Edges[style={dotted}](q3,q8) 
\Edges[style={dotted}](p2,p7)

\foreach \k in {p8,q2} {  \draw (\k) \kwb  ;}
\foreach \k in {p0,p2,p7,q10} {  \draw (\k) \kwc  ;}
\foreach \k in {q1,q9} {  \draw (\k) \kob  ;}
\foreach \k in {p1,p9,q0,q3,q8} {  \draw (\k) \koc  ;}

\draw (p0) node[below,yshift=-1pt]{$a_2$};
\draw (p1) node[below,yshift=-1pt]{$a_3$};
\draw (p2) node[below,yshift=-1pt]{$a_4$};

\draw (p8) node[below,yshift=-1pt]{$y_2$};
\draw (p9) node[below,yshift=-1pt]{$b_2$};

\draw (q0) node[above,yshift=1pt]{${a_1}$};
\draw (q1) node[above,yshift=1pt]{${x_1}$};
\draw (q2) node[above,yshift=1pt]{${x_2}$};
\draw (q3) node[above,yshift=1pt]{$a_5$};

\draw (q8) node[above,yshift=1pt]{$b_1$};
\draw (q9) node[above,yshift=1pt]{$y_1$};

\draw (q5) node[above,yshift=1pt]{$N_{b}$};
\draw (p5) node[above,xshift=-5mm]{$N_{a}$};

\end{tikzpicture}
\caption{Vertices $x_{1}$ and $y_1$ are coloured the same by $B_{1}$}
\end{figure}
} 

\def\nowyviii{
\begin{figure}[!htb]
\centering

\begin{tikzpicture}
\tikzset{VertexStyle/.style = {
shape = circle,
fill = black,
inner sep = 0pt,
outer sep = 0pt,
minimum size = 0pt,
draw}}

\SetVertexNoLabel

\grEmptyPath[Math,prefix=p,RA=1,RS=0]{10}
\begin{scope}[xshift=0.5 cm]
\grEmptyPath[Math,prefix=q,RA=1,RS=1]{11}
\end{scope}

\Edges(p0,q0,p1,q1,p2,q2,p3,q3)
\Edges(q0,q3)
\Edges(p0,p3)
\Edges(q6,p7,q7,p8,q8,p9,q9)
\Edges(q6,q9)
\Edges(p7,p9)
\Edges[style={dotted}](q3,q6)
\Edges[style={dotted}](p3,p7)

\foreach \k in {p1,q2,q7} {  \draw (\k) \kwb  ;}
\foreach \k in {p3,p8,q0,q6,q9} {  \draw (\k) \kwc  ;}
\foreach \k in {q1,q8} {  \draw (\k) \kob  ;}
\foreach \k in {p0,p2,p7,p9,q3} {  \draw (\k) \koc  ;}

\draw (p1) node[below,yshift=-1pt]{$x_1$};
\draw (p2) node[below,yshift=-1pt]{$b_1$};

\draw (p7) node[below,yshift=-1pt]{$a_2$};
\draw (p8) node[below,yshift=-1pt]{$a_3$};
\draw (p9) node[below,yshift=-1pt]{$a_4$};

\draw (q1) node[above,yshift=1pt]{$x_2$};
\draw (q2) node[above,yshift=1pt]{$x_3$};

\draw (q6) node[above,yshift=1pt]{$a_1$};
\draw (q7) node[above,yshift=1pt]{$y_1$};
\draw (q8) node[above,yshift=1pt]{$y_2$};
\draw (q9) node[above,yshift=1pt]{$a_5$};

\draw (q4) node[above,xshift=5mm]{$N_{b}$};
\draw (p5) node[above,xshift=0mm]{$N_{a}$};
\end{tikzpicture}
\caption{Vertices $x_2$ and $y_2$ are coloured the same by $B_{1}$}
\end{figure}
} 

\def\nowyix{
\begin{figure}[!htb]
\centering

\begin{tikzpicture}
\tikzset{VertexStyle/.style = {
shape = circle,
fill = black,
inner sep = 0pt,
outer sep = 0pt,
minimum size = 0pt,
draw}}

\SetVertexNoLabel

\grEmptyPath[Math,prefix=p,RA=1,RS=0]{10}
\begin{scope}[xshift=0.5 cm]
\grEmptyPath[Math,prefix=q,RA=1,RS=1]{11}
\end{scope}

\Edges(p0,q0,p1,q1,p2,q2,p3,q3)
\Edges(q0,q3)
\Edges(p0,p3)
\Edges(q7,p8)
\Edges[style={dotted}](q3,q7)
\Edges[style={dotted}](p3,p8)

\foreach \k in {p1,p8,q2} {  \draw (\k) \kwb  ;}
\foreach \k in {p3,q0} {  \draw (\k) \kwc  ;}
\foreach \k in {q1,q7} {  \draw (\k) \kob  ;}
\foreach \k in {p0,p2,q3} {  \draw (\k) \koc  ;}

\draw (p1) node[below,yshift=-2pt]{$x_1$};
\draw (p2) node[below,yshift=-1pt]{$b_1$};

\draw (p8) node[below,yshift=-1pt]{$y_2$};

\draw (q0) node[above,yshift=1pt]{$b_2$};
\draw (q1) node[above,yshift=1pt]{$x_2$};
\draw (q2) node[above,yshift=1pt]{$x_3$};

\draw (q7) node[above,yshift=1pt]{$y_1$};

\draw (q5) node[above,yshift=1pt]{$N_{b}$};
\draw (p6) node[above,xshift=-5mm]{$N_{a}$};
\end{tikzpicture}
\caption{A pair of edges $(x_{1}x_{2}, y_{1}y_{2})$ splits the vertex set of $A(3, 4)$ into two parts}
\end{figure}
} 

\def\nowyx{
\newbox\Bo
\savebox\Bo{
\begin{tikzpicture}
\tikzset{VertexStyle/.style = {
shape = circle,
fill = black,
inner sep = 0pt,
outer sep = 0pt,
minimum size = 0pt,
draw}}

\SetVertexNoLabel

\grPath[Math,prefix=p,RA=1,RS=0]{5}
\begin{scope}[xshift=-0.5 cm]
\grPath[Math,prefix=q,RA=1,RS=1]{5}
\end{scope}

\Edges(q0,p0,q1,p1,q2,p2,q3,p3,q4,p4)
\Edges(q0,q4)
\Edges(p0,p4)

\foreach \k in {p0,p2,q4} {  \draw (\k) \kwc  ;}
\foreach \k in {p1,p3,q0} {  \draw (\k) \kwb  ;}
\foreach \k in {q2} {  \draw (\k) \kob  ;}
\foreach \k in {p4,q1,q3} {  \draw (\k) \koc  ;}

\foreach \k in{0,1,2,3,4} {
\draw (q\k) node[above]{$a_{\the\numexpr\k+1\relax}$};
\draw (p\k) node[below]{$b_{\the\numexpr\k+1\relax}$};}

\end{tikzpicture}
}

\newbox\Bi
\savebox\Bi{
\begin{tikzpicture}
\tikzset{VertexStyle/.style = {
shape = circle,
fill = black,
inner sep = 0pt,
outer sep = 0pt,
minimum size = 0pt,
draw}}

\SetVertexNoLabel

\grPath[Math,prefix=p,RA=1,RS=0]{5}
\begin{scope}[xshift=-0.5 cm]
\grPath[Math,prefix=q,RA=1,RS=1]{5}
\end{scope}

\Edges(q0,p0,q1,p1,q2,p2,q3,p3,q4,p4)
\Edges(q0,q4)
\Edges(p0,p4)

\foreach \k in {p0,p2,q4} {  \draw (\k) \kwc  ;}
\foreach \k in {p1,p3,q0} {  \draw (\k) \kwb  ;}
\foreach \k in {q1,q3} {  \draw (\k) \kob  ;}
\foreach \k in {p4,q2} {  \draw (\k) \koc  ;}

\end{tikzpicture}
}

\newbox\Bii
\savebox\Bii{
\begin{tikzpicture}
\tikzset{VertexStyle/.style = {
shape = circle,
fill = black,
inner sep = 0pt,
outer sep = 0pt,
minimum size = 0pt,
draw}}

\SetVertexNoLabel

\grPath[Math,prefix=p,RA=1,RS=0]{5}
\begin{scope}[xshift=-0.5 cm]
\grPath[Math,prefix=q,RA=1,RS=1]{5}
\end{scope}

\Edges(q0,p0,q1,p1,q2,p2,q3,p3,q4,p4)
\Edges(q0,q4)
\Edges(p0,p4)

\foreach \k in {p0,q2,q4} {  \draw (\k) \kwc  ;}
\foreach \k in {p1,p3,q0} {  \draw (\k) \kwb  ;}
\foreach \k in {q1,q3} {  \draw (\k) \kob  ;}
\foreach \k in {p2,p4} {  \draw (\k) \koc  ;}

\end{tikzpicture}
}

\newbox\Biii
\savebox\Biii{
\begin{tikzpicture}
\tikzset{VertexStyle/.style = {
shape = circle,
fill = black,
inner sep = 0pt,
outer sep = 0pt,
minimum size = 0pt,
draw}}

\SetVertexNoLabel

\grPath[Math,prefix=p,RA=1,RS=0]{5}
\begin{scope}[xshift=-0.5 cm]
\grPath[Math,prefix=q,RA=1,RS=1]{5}
\end{scope}

\Edges(q0,p0,q1,p1,q2,p2,q3,p3,q4,p4)
\Edges(q0,q4)
\Edges(p0,p4)

\foreach \k in {p0,q2,q4} {  \draw (\k) \kwc  ;}
\foreach \k in {p1,p3,q0} {  \draw (\k) \kwb  ;}
\foreach \k in {q3} {  \draw (\k) \kob  ;}
\foreach \k in {p2,p4,q1} {  \draw (\k) \koc  ;}

\end{tikzpicture}
}

\newbox\Bnic
\savebox\Bnic{
\begin{tikzpicture}
\tikzset{VertexStyle/.style = {
shape = circle,
fill = black,
inner sep = 0pt,
outer sep = 0pt,
minimum size = 0pt,
draw}}

\SetVertexNoLabel

\grPath[Math,prefix=p,RA=1,RS=0]{5}
\begin{scope}[xshift=-0.5 cm]
\grPath[Math,prefix=q,RA=1,RS=1]{5}
\end{scope}

\Edges(q0,p0,q1,p1,q2,p2,q3,p3,q4,p4)
\Edges(q0,q4)
\Edges(p0,p4)

\foreach \k in {q0,p1,p3} {  \draw (\k) \koc  ;}
\foreach \k in {p0,q2,q4} {  \draw (\k) \kwc  ;}
\foreach \k in {q3} {  \draw (\k) \kob  ;}
\foreach \k in {q1,p2,p4} {  \draw (\k) \kwb  ;}

\end{tikzpicture}
}

\begin{figure}[!htbp]
\centering
\begin{tikzpicture}
\node (B0) {\usebox\Bo};
\node[below=of B0] (B1) {\usebox\Bi};
\node[below=of B1] (B2) {\usebox\Bii};
\node[below=of B2] (B3) {\usebox\Biii};
\node[below=of B3] (Bnic) {\usebox\Bnic};
\draw[->] (B0) -- (B1) node[pos=.5,right] {$B_1$};
\draw[->] (B1) -- (B2) node[pos=.5,right] {$B_2$};
\draw[->] (B2) -- (B3) node[pos=.5,right] {$B_3$};
\draw[->] (B3) -- (Bnic) node[pos=.5,right] {$B$};

\end{tikzpicture}
\caption{$e = a_{3}b_{2}$ and $f = a_{4}b_{3}$  are consecutive parallel edges of type $1$}
\end{figure}
} 

\begin{document}
\begin{frontmatter}





\title{Kempe equivalence of $4$-colourings of some  plane triangulations} 
\author{Jan~Florek}
\ead{jan.florek@pwr.edu.pl}

\address{Faculty of Pure and Applied Mathematics,

 Wroclaw University of Science and Technology,

 Wybrze\.{z}e Wyspia\'nskiego 27,
50-370 Wroc{\l}aw, Poland}

\fntext[]{This research did not receive any specific grant from funding agencies in the public, commercial, or not-for-profit sectors.}

\begin{abstract}

Let $G_{n}$, where $n \geqslant 5$, be a simple plane triangulation which has $2$ non-adjacent vertices of degree $n$ (called  \textit{poles} of $G_n$) and $2n$ vertices of degree~$5$. 
A set of Kempe equivalent $4$-colourings of $G_{n}$ is called a \textit{Kempe class}. The number of Kempe classes of $G_{n}$ is enumerated. In particular it is shown that there is at least $\lfloor \frac{n}{6} \rfloor$ Kempe classes of $G_{n}$. 

We say that $4$-colourings $A, B$ of $G_{n}$ are \textit{equal} if there exists a permutation~$P$ of the set of colours  such that $A = P \circ B$. Otherwise, $A$, $B$ are  \textit{different}. The number of different $4$-colourings of  $G_{n}$ is enumerated.

Suppose that $H_{n} = G_{n} - b$, where $b$ is a pole of $G_{n}$. We prove that all $4$-colourings of  $H_{n}$ are Kempe equivalent up to  
 $\lfloor \frac{13n}{2} \rfloor$ Kempe changes.
\end{abstract}


\end{frontmatter}
\section{Introduction}
We use Bondy and Murty  \cite{flobar4} as a reference for undefined terms.

Let $G$ be a graph and $k \geqslant 1$ be an integer. A vertex set $U \subseteq V(G)$ is \textit{independent} if no two vertices are adjacent in $G$. A $k$-colouring of $G$ is a partition of $V(G)$ into $k$ independent sets $U_{1}$, \ldots, $U_{k}$ called \textit{colour classes}. If $v \in U_{i}$ ($i = 1,$ \ldots, $k$), then $v$ is said to have \textit{colour} $i$. Every $k$-colouring can be identified with a function $A \colon  V(G) \rightarrow \{1, \ldots, k\}$ such that $A(v)$ is the colour of $v$. For a colouring $A$ and distinct colours $i$, $j \in  \{1, \ldots, k\}$, $A_{G}(i, j)$ (shortly $A(i, j)$) is the subgraph of $G$ induced by vertices with colour $i$ and $j$. A component of $A(i, j)$ is called a \textit{Kempe chain}. A \textit{Kempe change} consists in swapping the two colours in a Kempe chain, thereby obtaining a new $4$-colouring of the graph.  A pair of colourings (say $A_1$ and~$A_2$) are \textit{Kempe equivalent} (in symbols $A_{1} \sim A_{2}$)  if one can be obtained from the other through a series of Kempe changes. Let ${\cal C}_{k}(G)$ be the set of all $k$-colourings of  $G$. A set of Kempe equivalent colourings of ${\cal C}_{k}(G)$ is called a \textit{Kempe class}. 

The study of Kempe changes has a vast history, see e.g. \cite{flobar15} and \cite{flobar2}. We briefly review studies of  Kempe equivalence. Fisk \cite{flobar10} showed that the set of all $4$-colourings of an Eulerian triangulation of the plane is a Kempe class. This was generalized both by Meyniel \cite{flobar13}, who showed that all $5$-colourings of a plane graph are Kempe equivalent, and by Mohar \cite{flobar15}, who proved that all  $k$-colourings of a plane graph $G$ are Kempe equivalent if $k > \chi(G)$, where $\chi(G)$ is the chromatic number of $G$. 
Las Vergnas and Meyniel \cite{flobar12}, showed that all $k$-colourings of a $d$-degenerate graph are equivalent for $k \geqslant d + 1$ (a~graph $G$, every subgraph of which has minimum degree at most $d$, is said to be \textit{$d$-degenerate}). Mohar \cite{flobar15} conjectured that all  $k$-colourings of a graph are Kempe equivalent for $k \geqslant \Delta$. Note that the result of Las Vergnas and Meyniel  settles the case of non-regular connected graphs.  Van den Heuvel \cite{flobar16} showed that there is a counterexample to the conjecture: the $3$-prism. Feghali \textit{et al.} \cite{flobar8} proved that the conjecture holds for all cubic graphs except of the $3$-prism. Bonamy  \textit{et al.} \cite{flobar2} affirmed the conjecture for $\Delta$-regular graphs with $\Delta > 4$.  Bonamy \textit{et al.} \cite{flobar3} proved that all $k$-colourings of an $n$-vertex graph $G$ with $\Delta \leqslant k$ are equivalent up to $O(n^{2})$ Kempe changes, unless $k = 3$ and $G$ is the $3$-prism. 
Deschamps et \textit{et al.} \cite{flobar7} proved that all $5$-colourings of an $n$-vertex plane graph are Kempe equivalent up to $O(n^{195})$ Kempe changes.

A $5$-connected plane triangulation is called \textit{essentially $6$-connected} if every separating $5$-cycle is induced by the set of neighbours of a vertex of degree~$5$ (see 
Bondy and Murty \cite{flobar2}). 
Let $G_{n}$, $n \geqslant 5$, be a simple plane triangulation which has two non-adjacent vertices of degree $n$ (called  \textit{poles} of $G_n$) and $2n$ vertices of degree $5$.  
Florek \cite{flobar11} proved that $\{G_{n}\colon n \geqslant 5\}$ is the family of all \textit{minimal} essentially $6$-connected triangulations which are not essentially $6$-connected  as soon as we contract an edge with an end-vertex of degree $5$. 

Fix $G_{n}$, for some $n \geqslant 5$. We say that colourings $A$, $B \in {\cal C}_{4}(G_{n})$ are \textit{equal} if there exists a permutation $P$ of the set $ \{1, 2, 3, 4\}$ such that $A = P \circ B$. Otherwise, $A$, $B$ are  \textit{different}. If $n \equiv 0\, (mod\, 3)$, then there exists exactly one $4$-colouring of $G_{n}$ (denoted by $Q$) which has both poles coloured the same. We may assume that poles of $G_{n}$ are coloured $1$ by $Q$. For every $A \neq Q$ we may assume that poles of $G_{n}$ are coloured $1$ and $2$. 

We say that an edge in $G_n$ is of \textit{type $1$} (of \textit{type $2$}) if its end-vertices are neighbours of different poles (of the same pole, respectively) of $G_{n}$. For every colouring $A \in {\cal C}_{4}(G_{n})$ we assign four numbers (see Definition \ref{definition1.1})  Namely, $a(A)$ (or $b(A)$) is the number of vertices of $V(G_n)$ coloured $1$ ($3$, respectively) by~$A$. $c(A)$ (or $d(A)$) is the number of edges of  type $2$ ($1$, respectively) in the subgraph $A(3, 4)$ ($A(1, 2)$, respectively). Moreover, we put $a(Q) = \frac{n}{3} +1$, $b(Q) = c(Q) = \frac{2n}{3}$ and $d(Q) = 0$, for $n \equiv 0 \pmod3$.
A colouring of ${\cal C}_{4}(G_{n})$ is  \textit{constant} if it is not equivalent to any other $4$-colouring of $G_{n}$. In Theorem~\ref{theorem1.1} we prove that if $A$, $B$ are not constant, then $A \sim B$ if and only if $a(A) = a(B)$ ($b(A) = b(B)$, $c(A) = c(B)$ and $d(A) = d(B)$, respectively).
Moreover, $A \sim Q$ if and only if $d(A) = 0$. $A$ is constant if and only if $d(A) = 1$. It follows that the above four numbers are  invariant under the Kempe changes. 

Let  $K^{\star}(G_{n}, 4)$ be the number of Kempe classes of $G_{n}$, where $\star$ means that the set of all constant colourings of ${\cal C}_{4}(G_{n})$ is treated as one Kempe class.  
In Theorem  \ref{theorem1.2} it is proved that
 \[
\begin{array}{ll}
K^{\star}(G_{n}, 4) =
\left\{
\begin{array}{ll}
 \left\lfloor \dfrac{n}{6} \right\rfloor + 1& \quad \hbox{for } n \not \equiv 1 \pmod6,
\\[9pt]
\left\lfloor \dfrac{n}{6} \right\rfloor & \quad \hbox{for } n  \equiv 1
\pmod6.
\end{array}
\right.
\end{array}
\]
If  $n \equiv 2 \pmod3$,
then there exist $2n$ colourings of ${\cal C}_{4}(G_{n})$ which are constant (see condition $(b)$ of Lemma \ref{lemma1.4} and Remark~\ref{remark1.2}). In Theorem \ref{theorem1.3} the order of the family  ${\cal C}_{4}(G_{n})$ is enumerated.

In chapter $3$ we consider a graph $H_{n} = G_{n} - b$ where $b$ is a pole of $G_{n}$, for $n \geqslant 5$. In Theorem \ref{theorem2.1} we show that for every graph $H_{n}$, every two $4$-colourings of $H_{n}$ are equivalent up to 
\begin{align*}
&  6\left\lfloor \frac{n}{2} \right\rfloor \ \hbox{Kempe changes}, \ \hbox{for} \  n \equiv 0
\pmod3,
\\
& 9\left\lfloor \frac{n}{2} \right\rfloor \ \hbox{Kempe changes},  \ \hbox{for} \  n \equiv 2
\pmod3,
\\
&9\left\lfloor {\frac{n}{2}}\right\rfloor + 6\left\lfloor {\frac{n}{3}}\right\rfloor - 2\ \hbox{Kempe changes},  \ \hbox{for} \  n \equiv 1
\pmod3.
\end{align*}

\section{Kempe invariants and Kempe equivalence classes of $G_{n}$}
Fix $G_{n}$, for some $n \geqslant 5$. Let $a$, $b$ be poles of $G_{n}$. Recall that ${\cal C}_{4}(G_{n})$ is the set of all $4$-colourings of the graph $G_{n}$. We assume that both poles of $G_{n}$ are coloured $1$ by $Q$. For every $A \in {\cal C}_{4}(G_{n})$, $A \neq Q$, we assume that poles are coloured $1$ and $2$ by $A$. 

For $A \in {\cal C}_{4}(G_{n})$, if $v$ is a vertex of $G_{n}$ indicated in Figs 1, 2 by white circle (white square, black circle and black square), then $A(v) = 1$ ($A(v) = 2$, $A(v) = 3$ and $A(v) = 4$, respectively). 

Let $N_{a}$ ($N_{b}$) be a clockwise oriented cycle induced by all neighbours of the pole $a$ ($b$, respectively). We may assume that $a$  belongs to the bounded region of $R^{2}\setminus N_{a}$.

 \begin{definition}\label{definition1.1} 
Let $A \in {\cal C}_{4}(G_{n})$. We say that an edge in $G_n$ is of \textit{type $1$} (of \textit{type $2$}) if one of its vertices belongs to $N_a$ and the other to $N_{b}$ (the edge is contained in $N_{a} \cup N_{b}$, respectively). A path or cycle in $G_{n}$ is called of \textit{kind $1$} (of \textit{kind $2$}) if its edges are of type $1$ and type $2$ alternately (its all edges are of type $2$, respectively). 
\end{definition}

\begin{lemma}\label{lemma1.1} Let $A \in {\cal C}_{4}(G_{n})$. Suppose that $\xi$ is a Kempe chain of the colouring $A$ not containing any pole of $G_n$. If $\xi$ is a component of $A(3,4)$,  then it is a path or a cycle of kind $1$ of even order. Moreover,  if  it is a path, then it is an edge of type $1$ or its both end-edges are of type $1$. If $\xi$ contains a vertex coloured $1$ or $2$, then it is a path or a cycle of kind $2$ of even order.
\end{lemma}

\begin{PrfFact} Let $A \in {\cal C}_{4}(G_{n})$ and suppose that $\xi$ is a Kempe chain of the colouring $A$ not containing any pole of $G_n$. 

If $\xi$ is a component of $A(3,4)$,  then it does not contain two adjacent edges of type $2$, because $A(1, 2)$ contains no edge of type $2$. Hence, $\xi$ is a path or a cycle of  kind $1$. Similarly,  if  it is a path, then it is an edge of type $1$ or its both end-edges  are of type $1$. Hence, $\xi$ is of even order.

If  $\xi$ is a component of $A(1, 3)$, then it is a path or a cycle of kind~$2$, because $\xi$ contains no pole of $G_n$. If it is a path, then the set of all neighbours of the vertex set of  $\xi$ induces a cycle contained in $A(2,4)$ which has four vertices more than $\xi$. Hence, $\xi$ is of even order. Similarly, if $\xi$ is a component of $A(1, 4)$ ($A(2, 3)$ and $A(2, 4)$, then $\xi$ is a path or a cycle of kind $2$ of even order.
\end{PrfFact}

Since $N_a$ ($N_{b}$) has the clockwise orientation, we may enumerate consecutive neighbours of $a$ ($b$, respectively), consecutive edges of type $1$ or type $2$ around the pole $a$ ($b$, respectively).
\begin{lemma}\label{lemma1.2}  Let $A \in {\cal C}_{4}(G_{n})$. Then, 
the numbers of vertices in $G_n$ coloured $1$ and $2$ ($3$ and $4$) are equal. 
\end{lemma}
\begin{PrfFact} 
Notice that, by Lemma \ref{lemma1.1}, each component of $A(3, 4)$ is a path or a cycle of kind $1$ of even order. Hence,  the numbers of vertices in $G_n$ coloured $3$ and $4$ are equal. 

If ${\xi}_1$ and ${\xi}_2$ are two consecutive components in $A(3, 4)$ each of order at least $4$, then, by Lemma \ref{lemma1.1}, the last edge of ${\xi}_1$ and the first edge of ${\xi}_2$ are of type $1$.  Hence, the last edge of type $2$ in ${\xi}_1$ belongs to $N_a$ if and only if  the first edge of type $2$ in ${\xi}_2$ belongs to $N_b$.  Thus, consecutive  edges  of type~$2$ in $A(3, 4)$  belong to $N_a$ and $N_b$ alternately. Hence, we obtain
\begin{enumerate} 
\item[$(i)$] the numbers of edges of type $2$ in $A(3, 4) \cap N_a$ and $A(3, 4) \cap N_b$ are equal
\end{enumerate}
Certainly, we may assume that the pole $a$ is coloured $1$ and $b$ is coloured $2$. Then, a vertex of $N_a$ is coloured $2$ if and only if it is a vertex of some edge of type $1$ in $A(1, 2)$, or it is adjacent to both end-vertices of some edge of type $2$ in  $A(3, 4) \cap N_b$.  Similarly, a vertex of $N_b$ is coloured $1$ if and only if it is a vertex of some edge of type $1$ in $A(1, 2)$, or it is adjacent to both end-vertices of some edge of type $2$ in  $A(3, 4) \cap N_a$.  Hence, by condition $(i)$, the numbers of vertices in $G_n$ coloured $1$ and $2$ are equal. 
\end{PrfFact}

\begin{definition}\label{definition1.2}  Let $A \in {\cal C}_{4}(G_{n})$, $A \neq Q$. $a(A)$ (or $b(A)$) denotes the number of vertices of $V(G_n)$ coloured $1$ ($3$, respectively) by $A$. Moreover, $c(A)$ (or $d(A)$) denotes the number of edges of  type $2$ ($1$, respectively) in $A(3, 4)$  ($A(1, 2)$, respectively). Further, we put $a(Q) := \frac{n}{3} +1$, $b(Q) = c(Q) = \frac{2n}{3}$ and $d(Q) = 0$, for $n \equiv 0\, 
\pmod3$.
\end{definition}

\begin{lemma}\label{lemma1.3} Let $A \in {\cal C}_{4}(G_{n})$.    
The following equations are satisfied:
\begin{enumerate}
\item[$(1)$] $a(A) +b(A) = n+1$,
\item[$(2)$] $c(A) +d(A) = b(A)$,
\item[$(3)$] $c(A) + 2d(A) = 2a(A) - 2$,
\item[$(4)$] $3b(A) + d(A) = 3c(A) + 4d(A) = 2n$.
\end{enumerate}
\end{lemma}

\begin{PrfFact} Certainly, if $A = Q$ then lemma holds. 

Let now $A \neq Q$. By Lemma \ref{lemma1.2}, $A(1, 2)$ (or  $A(3, 4)$) has $2a(A)$ ($2b(A)$, respectively)  vertices. Hence, condition $(1)$ holds.

 By Lemma \ref{lemma1.1} each component of $A(3,4)$ is a path or cycle of kind $1$. Moreover,  if  it is a path, then its end-edges are both of type $1$. Hence, each vertex of $A(3, 4)$ satisfies exactly one of the following conditions: 
\begin{enumerate}
\item[$(i)$] it is a vertex of an edge of type $2$ in $A(3, 4)$, 
\item[$(ii)$] it is adjacent to both end-vertices of some edge of type $1$ in $A(1, 2)$.
\end{enumerate}
Since $A(3, 4)$ has $2b(A)$ vertices condition $(2)$ holds. 

Notice that each vertex of $A(1,2)$ different from a pole, satisfies exactly one of the following conditions: 
\begin{enumerate}
\item[$(iii)$] it is adjacent to both end-vertices of some edge of type $2$ in $A(3, 4)$,
\item[$(iv)$] it is a vertex of an edge of type $1$ in $A(1, 2)$.
\end{enumerate}
Since $A(1, 2)$) has $2a(A)$ vertices condition $(3)$ holds. 

By conditions $(2)$, $(3)$ and $(1)$, we obtain
\begin{multline*}
3b(A) + d(A) = 3(c(A) + d(A)) + d(A) = 3c(A) + 4d(A) =
\\ 
= \, 2c(A) + 2d(A) + c(A) + 2d(A) = 2b(A) + 2a(A) - 2 = 2n.
\end{multline*}
\end{PrfFact}

\begin{lemma}\label{lemma1.4} Let $A \in {\cal C}_{4}(G_{n})$ with $d(A) > 1$. 
If $B \in {\cal C}_{4}(G_{n})$ and $B \sim A$, then $d(B) = d(A)$.
\end{lemma}

\begin{PrfFact} Let $A \in {\cal C}_{4}(G_{n})$ with $d(A) > 1$ and suppose that $\xi$ is a proper Kempe chain contained in $A(i, j)$, for some different $i, j \in \{1, 2, 3, 4\}$.
Assume that $B$ is a colouring obtained from $A$ by switching two colours in $\xi$. By Lemma~\ref{lemma1.1} one of the following conditions is satisfied:
\begin{enumerate}
\item[$(i)$] $\xi$ is a path of kind $1$ in $A(3, 4)$,
\item[$(ii)$] $\xi$ is a path of kind $2$ of even order containing a vertex coloured $1$ or $2$, 
\item[$(iii)$] $\xi$ contains a pole of $G_{n}$.
\end{enumerate}
Case $(i)$. Then, $d(B) = d(A)$.

Case $(ii)$. Assume that $\xi$ is a path of kind $2$ in $A(2,4)$. Let $A_{\xi}(1, 2)$ (or $A_{\xi}(1, 4)$) be the set of edges of type $1$ in $A(1, 2)$ (or $A(1, 4)$) with one end-vertex belonging to $\xi$.  Since $\xi$ is of even order, $|A_{\xi}(1, 2)| = |A_{\xi}(1, 4)|$. Since $B$ is a colouring obtained from $A$ by switching colours $2$ and $4$ in $\xi$, then $A_{\xi}(1, 4) = B_{\xi}(1, 2)$.  Hence, $|A_{\xi}(1, 2)| = |A_{\xi}(1, 4)|= |B_{\xi}(1, 2)|$. Therefore, $d(A) = d(B)$.

Similarly, if  $\xi$ is a path of kind $2$ in $A(2,3)$ ($A(1,3)$ and  $A(1, 4)$) then $d(B) = d(A)$.

Case $(iii)$. Since $d(A) > 0$,  $A(1, 2)$ is connected. Thus  $\xi$ is not a proper Kempe chain in $A(1, 2)$.

Assume that $\xi$ is a proper Kempe chain in $A(2, 4)$ containing the pole coloured  $2$.  If $B' \in {\cal C}_{4}(G_{n})$  is a colouring obtained from $A$ by switching the colours in each component  of $A(2, 4)$ different from $\xi$, then $B'$ is equal to $B$. Notice that each component of $A(2, 4)$ different from $\xi$ does not  contain the pole coloured $2$. By Lemma \ref{lemma1.1}, each of them is a path of kind $2$ of even order. Hence,  $d(B') = d(A)$, by condition $(ii)$. Thus, $d(B) = d(B') = d(A)$.

Similarly, if  $\xi$ is a component of $A(2, 3)$ ($A(1, 3)$ and $A(1, 4)$) containing the pole coloured  $2$ (coloured $1$, respectively), then $d(B) = d(A)$.
 \end{PrfFact}
 
\begin{lemma}\label{lemma1.5} 
Let $A \in {\cal C}_{4}(G_{n})$.
 \begin{enumerate} 
  \item[$(a)$] $A$ is constant if and only if $d(A) = 1$,
   \item[$(b)$] $A \sim Q$ if and only if $d(A) = 0$. $\{A \in {\cal C}_{4}(G_{n}): A \sim Q\}$ has four elements.
 \end{enumerate}
  \end{lemma}
  
  \begin{PrfFact} 
$(a)$ Certainly, if $A$ is constant, then $d(A) = 1$. 

Let $d(A) = 1$. Then, $A(1, 2)$ and $A(3, 4)$ has no proper Kempe chain. If $A(2, 4)$ contains a proper Kempe chain (say $\xi$), then, by Lemma \ref{lemma1.1}, it is a path of kind $2$ of even order. If $\xi$ is of length at least $3$, then it contains at least $2$ vertices coloured $2$. Hence, $A(1, 2)$ contains at least $2$ edges of type~$1$ which is a contradiction. If $\xi$ is an edge, then $A(3, 4)$ is a path of kind $1$ of odd order  (see Fig. 2) which, by Lemma \ref{lemma1.1}, is a contradiction. Similarly, $A(2, 3)$ ($A(1, 4)$ and $A(1, 3)$) contains no proper Kempe chain. Hence, $A$ is constant and the condition $(a)$ holds.
 
Proof $(b)$ Let both poles of $G_{n}$ be coloured $1$ by $Q$. Then, $G_n$ has exactly three different cycles of kind $1$: $Q(3, 4)$, $ Q(4, 2)$ and $Q(2, 3)$. It follows that there are exactly three colourings of ${\cal C}_{4}(G_{n})$ (say $A, B, C$) different from $Q$ which can be obtained from $Q$ by a single Kempe change. 
  
Let now $D  \in {\cal C}_{4}(G_{n)}$, $D \neq Q$,  and suppose that poles of $G_{n}$ are coloured~$1$ and $2$. Assume that $D \sim Q$. By Lemma \ref{lemma1.4} and condition $(a)$, $d(D) = 0$. Hence, $D(3, 4)$ is a cycle of type $1$. Then, $D$  can be obtained from $Q$ by a single Kempe change. Therefore, $D \in \{A, B, C\}$ and the condition $(b)$ holds.
 \end{PrfFact} 
 
 \nowyi
 
 \begin{definition}\label{definition1.3} Let $e$ be an edge of type $1$ and $k$ be an integer, $1 \leqslant k <~\frac{n}{2}$. A colouring $A  \in {\cal C}_{4}(G_{n})$ is denoted by $Q_{k, e}$ if $d(A) = k$ and  $A(3, 4)$ has $k - 1$ consecutive components which are edges of type $1$ and $e$ is the first of these edges. A colouring $A  \in {\cal C}_{4}(G_{n})$ is denoted by  $Q_{1, e}$ (or  $Q_{\frac{n}{2}, e}$)  if $d(A) = 1$ ($d(A) =\frac{n}{2}$, respectively), $A(3, 4)$  is a path of type $1$ and $e$ is its first edge ($ A(3, 4)$ has $\frac{n}{2}$ components each of which is an edge of type $1$ and $e$ is one of them, respectively). Notice that $Q_{k, e}$ is not defined clearly (there exist $2^{k-1}$ different $4$-colourings of $G_{n}$ called  $Q_{k, e}$).
 \end{definition}
  
  \begin{remark}\label{remark1.2}  Let $e$ be an edge of type $1$ and $1 \leqslant k \leqslant~\frac{n}{2}$. If there exists a colouring $Q_{k, e}$ of  $G_{n}$, then, by condition $(4)$ of Lemma \ref{lemma1.3},  $n \equiv 2k \, 
\pmod3$.
It is easy to see that if $n  \equiv 2k \,
\pmod3$,
then there exists $Q_{k, e}$ (see Fig. 1).  
 \end{remark}
 
 \begin{lemma}\label{lemma1.6} 
 Let $A  \in {\cal C}_{4}(G_{n})$ with $d(A) = k > 1$.  
 Assume that $\xi_{1}$, $\xi_{2}$, \ldots, $\xi_{k}$ is a sequence  of $k$ consecutive components of $A(3,4)$ such that that $|\xi_{1}| \geqslant 4$. Then, there exists a colouring $A' \in {\cal C}_{4}(G_{n})$  equivalent to $A$ such that $A^{'}(3,4)$ has $k$ consecutive components $\xi^{'}_{1}, \xi^{'}_{2}, \ldots, \xi^{'}_{k}$ satisfying the following conditions: 
\begin{enumerate}
\item[$(1)$] $|\xi^{'}_{1}| = |\xi_{1}| - 2$,
\item[$(2)$] $|\xi^{'}_{2}| = |\xi_{2}| + 2$,
\item[$(3)$] $\xi^{'}_{i} = \xi_{i}$, for $i > 2$,
\item[$(4)$] $\xi^{'}_{1}$ and $\xi_{1}$ have the same first edge,
\item[$(5)$] the last edge of $\xi_{1}$ and the first edge of $\xi^{'}_{2}$ are consecutive edges of type~$1$ in $G_{n}$,
\item[$(6)$] the first edge of $\xi_{2}$ is the third edge of $\xi^{'}_{2}$.
\end{enumerate} 
 \end{lemma}
 
  \nowyii
  
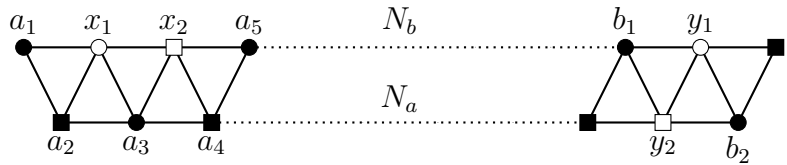
 \begin{PrfFact} Notice that, by Lemma \ref{lemma1.1}, $\xi_{i}$, for $i =1$, \ldots, $k$, is an edge or a path type $1$ with both end-edges of type $1$. Let $c_{1}$, \ldots, $c_{p}$ be consecutive vertices of $\xi_{1}$, where $p \geqslant 4$,  and suppose that $d_{1}$, \ldots, $d_{r}$ are consecutive vertices of~$\xi_{2}$.  Assume that $c$ is a common neighbour of the vertices $c_{p - 1}$, $c_{p}$ and~$d_{1}$. Switching colours on vertices of $\xi_{2}$ we obtain a $4$-colouring of $G_n$ equivalent to $A$ such that  $c_p$ and $d_1$ are coloured the same. Hence, we assume that  vertices $c_p$ and $d_1$ are coloured $3$ and $c$ is coloured $2$ by $A$ (see Fig. 2). Notice that the edge $c_{p-1}c$  is a component of $A(2, 4)$. Switching colours in $c_{p-1}c$ we obtain a colouring $A^{'} \in  {\cal C}_{4}(G_{n})$ equivalent to $A$. Components of $A^{'}(3, 4)$ satisfy the conditions $(1)-(6)$. Hence, lemma holds. 
 \end{PrfFact}
 
 \begin{corollary}\label{corollary1.1}
Let $A  \in {\cal C}_{4}(G_{n})$ with $d(A) = k > 1$. Assume that  $\xi_{1}$, $\xi_{2}$, \ldots, $\xi_{k}$ is a sequence  of $k$ consecutive components of $A(3, 4)$. Then there exists a colouring $B_{1} \in {\cal C}_{4}(G_{n})$ equivalent to $A$ such that $B_{1}(3, 4)$ has $k$ consecutive components $\sigma_{1}$, $\sigma_{2}$, \ldots, $\sigma_{k}$ satisfying the following conditions: 
\begin{enumerate}
\item[$(1)$] $|\sigma_{1}| = 2$
\item[$(2)$] $|\sigma_{2}| = |\xi_{1}| + |\xi_{2}| - 2$,
\item[$(3)$] $\sigma_{i} = \xi_{i}$, for $i > 2$,
\item[$(4)$] $\sigma_{1}$ is the first edge of $\xi_{1}$.
\end{enumerate} 
  \end{corollary}
  
  \begin{PrfFact} If $|\xi_{1}| = 2$ then corollary holds.  If $|\xi_{1}|  \geqslant 4$, then the colouring~$A^{'}$ (defined in Lemma \ref{lemma1.6}) is equivalent to $A$ and $A^{'}(3.4)$ consists of $k$ consecutive components  $\xi^{'}_{1}$, $\xi^{'}_{2}$, \ldots, $\xi^{'}_{k}$   satisfying conditions $(1) - (6)$ of Lemma \ref{lemma1.6}. If $|\xi^{'}_{1}| = 2$ then corollary holds. If $|\xi^{'}_{1}| \geqslant 4$ we continue the process. Finally, we obtain a colouring $B_{1} \in {\cal C}_{4}(G_{n})$ equivalent to $A$ such that $B_{1}(3, 4)$ has $k$  components satisfying the conditions (1)--(4). 
  \end{PrfFact}
  
 \begin{corollary}\label{corollary1.2} Let $A  \in {\cal C}_{4}(G_{n})$ with $d(A) = 2$.  Assume that $\xi_{1}$, $\xi_{2}$ are components of $A(3, 4)$  such that $|\xi_{1}| \leqslant |\xi_{2}|$, $e$ is the first edge of $\xi_{1}$ and  the last vertex of $\xi_{1}$ and the first vertex of $\xi_{2}$ are coloured the same by $A$. Then $A$ and $Q_{2, e}$ are equivalent up to $\frac{c(A)}{2}$ Kempe changes each of which switches colours in some edge of type~$2$.
 \end{corollary}
 
   \begin{PrfFact} If $|\xi_{1}| = 2$ then $A = Q_{2, e}$ and corollary holds.  If $|\xi_{1}| \geqslant 4$, then the colouring $A^{'}$ (defined in Lemma \ref{lemma1.6}) is obtained from $A$ by switching the colours in the edge $c_{p-1}c$ of type $2$. Notice that the edge $c_{p-1}c$ is incident with the last edge of type $2$ contained in $\xi_{1}$ (the edge $c_{p-2}c_{p-1}$). Then, $A^{'}(3, 4)$ consists of two components $\xi^{'}_{1}$, $\xi^{'}_{2}$ such that the last vertex of $\xi^{'}_{1}$  and the first vertex of $\xi^{'}_{2}$ are coloured the same by $A^{'}$ (vertices $c_{p - 2}$ and $c_{p}$). Moreover, $c(\xi^{'}_{1}) = c(\xi_{1})- 1$. If  $|\xi^{'}_{1}| = 2$, then $A^{'} = Q_{2, e}$. If $|\xi^{'}_{1}| \geqslant 4$ we continue the process. Finally, we obtain a colouring $Q_{2, e}$ after $c(\xi_{1}) \leqslant \frac{c(A)}{2}$ Kempe changes (where $c(\xi_{1})$ is the number of edges of type $2$ in $\xi_{1}$)  each of which switches colours in some edge of type $2$.
  \end{PrfFact}
  
 \begin{lemma}\label{lemma1.7}
Let $A  \in {\cal C}_{4}(G_{n})$ with $d(A) = k > 1$. Assume that  $\xi_{1}$, $\xi_{2}$, \ldots, $\xi_{k}$ is a sequence of $k$ consecutive components of $A(3, 4)$ such that  $|\xi_{1}| \geqslant 4$ and $e$ is the first edge of $\xi_{1}$. Then $A \sim  Q_{k, e}$. 
  \end{lemma}
  
 \begin{PrfFact}
Let $j$ be a maximal integer such that there exists a colouring $B_{j} \in {\cal C}_{4}(G_{n})$ such that $B_{j} \sim A$ and $B_{j}(3, 4)$ has $k$ consecutive components $\sigma_{1}$, $\sigma_{2}$, \ldots, $\sigma_{k}$ satisfying the following conditions: 
\begin{enumerate}
\item[$(1)$] $\sigma_{i}$ is an edge, for $1 \leqslant i \leqslant j$,
\item[$(2)$] $|\sigma_{j+1}| \geqslant 4$,
\item[$(3)$] $\sigma_{i} = \xi_{i}$ for $i > j+1$,
\item[$(4)$] $\sigma_{1} = e$.
\end{enumerate}
We will prove that $j = k - 1$. If $k = 2$, then by Corollary \ref{corollary1.1}, there exists a colouring $B_{1} \in {\cal C}_{4}(G_{n})$ satisfying conditions (1)--(4). 

Let $k \geqslant 3$ and suppose, on the contrary, that $j < k-1$. Then, $\sigma_{j+1}$, $\sigma_{j+2}$, \ldots, $\sigma_{k}$, $\sigma_{1}$, \ldots, $\sigma_{j}$ is a sequence of consecutive components of $B_{j}(3, 4)$. By condition $(2)$, $|\sigma_{j+1}| \geqslant 4$. Hence, by Corollary \ref{corollary1.1}, there exists $B_{j+1} \in~{\cal C}_{4}(G_{n})$ such that $B_{j+1} \sim B_{j}$ and $B_{j+1}(3, 4)$ has $k$ consecutive components $\delta_{j+1}$, $\delta_{j+2}$, \ldots, $\delta_{k}$, $\delta_{1}$,  \ldots, $\delta_{j}$ satisfying the following conditions:
\begin{enumerate}
\item[$(5)$] $|\delta_{j+1}| = 2$,
\item[$(6)$] $|\delta_{j+2}| \geqslant 4$,
\item[$(7)$] $\delta_{i} = \sigma_{i}$ for $i \neq {j+1}$ and $i \neq {j+2}$.
\end{enumerate} 
By conditions $(7)$ and $(3)$, $\delta_{i} = \sigma_{i} = \xi_{i}$, for $ i  > j+2$. Moreover, by conditions $(7)$ and $(1)$, $\delta_{i} = \sigma_{i}$ is an edge, for $1 \leqslant i \leqslant j$. By condition~$(5)$, $\delta_{j+1}$ is an edge. Further, by conditions  $(7)$ and $(4)$, $\delta_{1} = e$. Hence, we obtain
\begin{enumerate}
\item[$(8)$] $\delta_{i}$ is an edge of type $1$, for $1 \leqslant i \leqslant j+1$,
\item[$(9)$] $|\delta_{j+2}| \geqslant 4$,
\item[$(10)$] $\delta_{i} = \xi_{i}$ for $i > j+2$,
\item[$(11)$] $\delta_{1}  = e$,
\end{enumerate}
which contradicts the maximality of $j$. Hence, $j = k - 1$ and, by condition~$(4)$, $\sigma_{1} = e$. Therefore, $A \sim B_{k-1}=Q_{k, e}$. 
\end{PrfFact}

 \begin{lemma}\label{lemma1.8} $Q_{k, e} \sim Q_{k, f}$ for every $1 < k < \frac{n}{2}$. 
 \end{lemma}
 
 \begin{PrfFact}
It is sufficient to prove the lemma when $e$ and $f$ are consecutive edges of type $1$ in $G_{n}$ (having a common vertex). Suppose that $\xi_{1}$, $\xi_{2}$, \ldots, $\xi_{k}$ is a sequence of consecutive components of $Q_{k, e}(3, 4)$ such that  $|\xi_{1}| \geqslant 4$, $\xi_{2} = e$ and $\xi_{i}$ is an edge, for $i > 1$. Hence, by Lemma \ref{lemma1.6}, there exits a colouring $Q' \in {\cal C}_{4}(G_{n})$ such that $Q' \sim Q_{k, e}$ and $Q'(3, 4)$ has $k$ consecutive components $\xi^{'}_{1}$, $\xi^{'}_{2}$, \ldots, $\xi^{'}_{k}$ satisfying the following conditions:
\begin{enumerate}
\item[$(1)$] $|\xi^{'}_{1}| = |\xi_{1}| - 2$,
\item[$(2)$] $|\xi^{'}_{2}| = |\xi_{2}| + 2$,
\item[$(3)$] $e$ is the third edge of  $\xi^{'}_{2}$.
\end{enumerate} 
Notice that  $\xi^{'}_{2}$, $\xi^{'}_{3}$, \ldots, $\xi^{'}_{k}$, $\xi^{'}_{1}$ ($\xi^{'}_{2}$, $\xi^{'}_{1}$, for $k = 2$) is a sequence of $k$ consecutive components of $Q'(3, 4)$. Since $\xi_{2} = e$, by condition $(2)$, $|\xi^{'}_{2}|= 4$. Hence, by  Lemma \ref{lemma1.6}, there exits a colouring $Q'' \in {\cal C}_{4}(G_{n})$ such that $Q'' \sim Q'$ and $Q''(3, 4)$ has a sequence of $k$ consecutive components $\xi^{''}_{2}$, $\xi^{''}_{3}$, \ldots, $\xi^{''}_{k}$, $\xi^{''}_{1}$ ($\xi^{''}_{2}$, $\xi^{''}_{1}$, for $k = 2$) satisfying the following conditions:
\begin{enumerate}
\item[$(4)$] $|\xi^{''}_{2}| = |\xi{'}_{2}| - 2$,
\item[$(5)$] $|\xi^{''}_{3}| = |\xi^{'}_{3}| + 2$,
\item[$(6)$] the last edge of $\xi^{'}_{2}$ and the first edge of $\xi^{''}_{3}$ are consecutive edges of type~$1$. 
\end{enumerate} 
Notice that $\xi^{''}_{3}$, \ldots, $\xi^{'''}_{k}$, $\xi^{''}_{1}$, $\xi^{''}_{2}$ ($\xi^{''}_{1}$, $\xi^{''}_{2}$, for $k = 2$) is a sequence of $k$ consecutive components of $Q^{''}(3, 4)$ and, by condition $(5)$, $|\xi^{''}_{3}| \geqslant 4$. Since $|\xi^{'}_{2}|= 4$, by condition $(3)$, $e$ is the last edge of $\xi^{'}_{2}$. Thus, by condition $(6)$, $f$ is the first edge of  $\xi^{''}_{3}$. Hence, by Lemma \ref{lemma1.7}, $Q'' \sim Q_{k, f}$ and $Q_{k, e} \sim Q' \sim Q'' \sim Q_{k, f}$.
\end{PrfFact}

  \begin{lemma}\label{lemma1.9} Let $A, B  \in {\cal C}_{4}(G_{n})$.
If $d(A) = d(B) > 1$, then $A \sim B$.
 \end{lemma}
 
 \begin{PrfFact}
Let $1 < d(A) = d(B) < \frac{n}{2}$. Then, $A(3, 4)$ ($B(3, 4)$) has $d(A) > 1$ components which are paths of kind $1$ and one of them contains at least $3$ edges. Let $e$ (or $f$, respectively) be the first edge of this component. Since  colourings $A$, $B$ satisfy assumption of  Lemma \ref{lemma1.7} we obtain $A \sim  Q_{d(A), e}$ and  $B \sim Q_{d(B), f}$. Hence, by Lemma \ref{lemma1.8},  $A \sim  Q_{d(A), e} \sim Q_{d(B), f} \sim B$. 
   
If $d(A) = d(B) = \frac{n}{2}$, then each Kempe chain of $A(3, 4)$ and $B(3, 4)$ is an edge of type $1$. We may switch colours on some edges of $A(3, 4)$ (or $B(3, 4)$) to obtain a colouring $A' \in {\cal C}_{4}(G_{n})$ ($B' \in {\cal C}_{4}(G_{n})$) such that $A'(2, 3) = N_{a}$ and $A'(1,4) = N_{b}$ ($B'(2, 3) = N_{a}$ and $B'(1,4) = N_{b}$, respectively). Certainly, $A' \sim B'$. Hence, $A \sim B$.  
 \end{PrfFact}
 
\begin{theorem}\label{theorem1.1}
For every two colourings $A$, $B  \in {\cal C}_{4}(G_{n})$ which are not constant the following conditions are equivalent: 
\begin{enumerate}
\item[$(1)$] $A \sim B$,
\item[$(2)$] $a(A) = a(B)$,
\item[$(3)$] $b(A) = b(B)$,
\item[$(4)$] $c(A) = c(B)$,
\item[$(5)$] $d(A) = d(B)$. 
\end{enumerate}
Moreover, if $A$, $B  \in {\cal C}_{4}(G_{n})$  are  constant, then conditions (2)--(5) are equivalent.
\end{theorem}
\begin{PrfFact} Let $A, B  \in {\cal C}_{4}(G_{n})$. If $d(A) > 1$, then, by Lemmas \ref{lemma1.4} and \ref{lemma1.9},  $A \sim B$ if and only if $d(A) = d(B)$. Notice that by Lemma \ref{lemma1.5}, $d(A) = 1$ if and only if $A$ is constant ($d(A) = 0$ if and only if $A \sim Q$). Hence, by Lemma~\ref{lemma1.3}, the theorem holds.
\end{PrfFact}

\begin{theorem}\label{theorem1.2} Let $K^{\star}(G_{n}, 4)$ denote the number of Kempe classes of $G_{n}$, where $\star$ means that the set of all constant colourings of ${\cal C}_{4}(G_{n})$ is treated as one Kempe class. 
 \[
 K^{\star}(G_{n}, 4) = \left|E\left[\frac{n}{2}, \frac{2n}{3}\right]\right| =
\left\{
\begin{array}{ll}
 \left\lfloor \dfrac{n}{6} \right\rfloor + 1& \quad \hbox{for } n \not \equiv 1\, 
 \pmod6,
\\[9pt]
\left\lfloor \dfrac{n}{6} \right\rfloor & \quad \hbox{for } n  \equiv 1\,
\pmod6,
\end{array}
\right.
\]
 where  $E[\frac{n}{2}, \frac{2n}{3}]$ is the set of all integers in the interval $[\frac{n}{2}, \frac{2n}{3}]$.
\end{theorem}

\begin{PrfFact}
We first prove that
\begin{enumerate}
\item[$(i)$] 
 a function $b\colon {\cal C}_{4}(G_{n}) \rightarrow E[\frac{n}{2}, \frac{2n}{3}]: A \rightarrow b(A)$ is a surjection.
\end{enumerate} 
According to conditions $(4)$ and $(2)$ of Lemma \ref{lemma1.3} we have 
\[
3b(A) \leqslant 3c(A) + 4d(A) = 2n \leqslant 4b(A).
\]
Hence, $b(A) \in E[\frac{n}{2}, \frac{2n}{3}]$. Let $l \in E[\frac{n}{2}, \frac{2n}{3}]$. 

If $\frac{2n}{3}$ is an integer, then $Q \in  {\cal C}_{4}(G_{n})$ and $b(Q) = \frac{2n}{3}$. 

If $\frac{n}{2} \leqslant l < \frac{2n}{3}$, then $1 \leqslant 2n - 3l \leqslant \frac{n}{2}$. Hence, by Remark \ref{remark1.2} there exists a colouring $Q_{2n-3l} \in  {\cal C}_{4}(G_{n})$. From condition $(4)$ of Lemma \ref{lemma1.3} we have
\[
b(Q_{2n-3l}) = \frac{2n -d(Q_{2n-3l})}{3} = \frac{2n - (2n - 3l)}{3} = l,
\]
 which yields $(i)$. Hence, by Theorem \ref{theorem1.1},
$ K^{\star}(G_{n}, 4) = |E[\frac{n}{2}, \frac{2n}{3}]|$.

It is easy to check that 
 \[
\left|E\left[\frac{n}{2}, \frac{2n}{3}\right]\right| =
\left\{
\begin{array}{ll}
\left \lfloor \dfrac{n}{6} \right\rfloor + 1& \quad \hbox{for } n \not \equiv 1\,
\pmod6,
\\[9pt]
\left\lfloor \dfrac{n}{6} \right\rfloor & \quad \hbox{for } n  \equiv 1\,
\pmod6,
\end{array}
\right.
\]
which completes the proof. 
\end{PrfFact}

 \begin{theorem} \label{theorem1.3}
For every $n \geqslant 3$ we have
\begin{align}
 |{\cal C}_{4}(G_{n})| =  \sum\limits_{k \in E[\frac{n}{2}, \frac{2n}{3})}  \binom{k}{2n - 3k}\frac{n2^{2n - 3k}}{k}, \hbox{ for} \quad  n \not \equiv 0 \, 
 \pmod3
\end{align}
 and \begin{align}
|{\cal C}_{4}(G_{n})| = \sum\limits_{k \in E[\frac{n}{2}, \frac{2n}{3})}  \binom{k}{2n - 3k}\frac{n2^{2n - 3k}}{k} +4, \hbox{ for} \quad  n \equiv 0 \, 
\pmod3, 
\end{align}
 where  $E[\frac{n}{2}, \frac{2n}{3})$ is the set of all integers in the interval $[\frac{n}{2}, \frac{2n}{3})$.
 \end{theorem}
 
\begin{PrfFact}
Let $A \in {\cal C}_{4}(G_{n})$, $A \neq Q$. If $A$ is not a constant colouring, then $[A]$ denotes the set of all colourings $B \in {\cal C}_{4}(G_{n})$ such that $B \sim A$. If $A$ is a constant colouring, then $[A]$ denotes the set of all constant colourings in ${\cal C}_{4}(G_{n})$. We first prove that 
\begin{align}
 \quad  \bigl| [A] \bigr| = \binom{c(A) + d(A) - 1}{d(A) - 1}\frac{n2^{d(A)}}{d(A)}.
 \end{align}
If $B \in [A]$, then by Theorem \ref{theorem1.1}, $d(B) = d(A)$ and $c(B) = c(A)$. Hence, we obtain 
\begin{enumerate}
\item[$(i)$] $B(3, 4)$ has $d(A)$ Kempe chains which are paths of kind $1$,
\item[$(ii)$] $B(3, 4)$ has $c(A)$ edges of type $2$.
\end{enumerate}
Fix an edge of type $1$ in $G_{n}$ (say $e$) and suppose that $[A, e]$ is the set of all colourings $B \in [A]$ such that $e$ is the first edge of some Kempe chain in $B(3, 4)$. By conditions $(i)$ and $(ii)$, 
\[
[A, e] \ \hbox{ has }\ S(A) 2^{d(A) -1} \ \hbox{ elements, where }  \ S(A) = \binom{c(A) + d(A) - 1}{d(A) - 1} 
\]
is the number of solutions in non-negative integers of the following equation 
$$x_{1} + \ldots + x_{d(A)}  = c(A).$$
Since $G_{n}$ has $2n$ edges and $e$ can   be the first edge of any of Kempe chain in $B(3,4)$, then equation $(3)$ holds.  

In view of condition $(2)$ and $(4)$ of Lemma \ref{lemma1.3}, $c(A) +d(A) = b(A)$ and  $3b(A) + d(A) = 2n$. Hence, by equation $(3)$, we obtain 
\begin{align*}
| [A] | =  \binom{b(A) - 1}{2n - 3b(A) - 1}\frac{n2^{2n - 3b(A)}}{2n - 3b(A)}
=  \binom{b(A)}{2n - 3b(A)}\frac{n2^{2n - 3b(A)}}{b(A)}.
\end{align*}
Thus, by Theorem \ref{theorem1.2}, equation $(1)$ holds for $n \not \equiv 0 \, 
\pmod3$.

If $Q \in  {\cal C}_{4}(G_{n})$, then,  by Lemma \ref{lemma1.5}$(a)$,  $[Q]$ has $4$ elements. Hence, equation $(2)$ holds for $n \equiv 0 \, 
\pmod3$.
\end{PrfFact}
\section{Kempe equivalence classes of $H_{n} = G_{n} - b$}
Fix $G_{n}$,  for some $n \geqslant 5$. Let $a$, $b$ be poles of $G_{n}$.  We recall that $N_{a}$ ($N_{b}$) is a clockwise oriented cycle induced by all neighbours of the pole $a$ ($b$, respectively). We may assume that $a$ belongs to the bounded region of $R^{2}\setminus N_{a}$.

Let $H_{n} = G_{n}- b$ be a subgraph of $G_{n}$. For every colouring $A \in {\cal C}_{4}(H_{n})$ we assume that $A(a) = 1$.  For distinct colours $i, j \in  \{1, 2, 3, 4\}$, $A_H(i, j)$ (shortly $A(i, j)$) is the subgraph of $H_{n}$ induced by vertices with colour $i$ and $j$.  The components  of $A(i, j)$ are called \textit{Kempe chains}. Each proper component of $A_H(i, j)$  is called a \textit{proper Kempe chain}. We say that colourings $A$, $B$ of ${\cal C}_{4}(H_{n})$ are \textit{Kempe equivalent} (in symbols $A \sim B$) if we can form one from the other by a sequence of Kempe changes. A Kempe change swapping two colours in a proper Kempe chain is called a \textit{proper Kempe change}. Notice that if a colouring $B \in {\cal C}_{4}(H_{n})$ is obtained from a colouring $A \in {\cal C}_{4}(H_{n})$ through a sequence of Kempe changes containing a subsequence of $m$ proper Kempe changes, then there exists a colouring  $B'$ equal to $B$ which is obtained from the colouring $A$ through a sequence of $m$ proper Kempe changes. Hence, if we bound the length of the shortest sequence of Kempe changes between any two colourings we may calculate only the number of proper Kempe changes of this sequence. 

For $A \in {\cal C}_{4}(H_{n})$, if $v$ is a vertex of $H_{n}$ indicated  in Figs~3, \ldots, 12 by white circle (white square, black circle and  black square) then $A(v) = 1$ ($A(v) = 2$,  $A(v) = 3$ and  $A(v) = 4$, respectively). 

Let $A \in {\cal C}_{4}(H_{n})$ and suppose $a_{1}b_{1}$, $a_{2}b_{2}$ are any disjoint edges such that $a_{1}, a_{2} \in N_a$ and $b_{1}, b_{2} \in N_b$. Then $C = aa_{1}b_{1}bb_{2}a_{2}a$ is a cycle in the graph $G_n$. If a subgraph $A(j, k)$ is disjoint with $C$, then we say that the pair $(a_{1}b_{1}, a_{2}b_{2})$ \textit{splits} the set $A(j, k)$  \textit{into two parts}: one part of  $A(j, k)$ is contained in the bounded and another one is contained in the unbounded region  determined by $C$ on the plane. 

We say that an edge in $H_n$ is of \textit{type $1$} (\textit{type $2$}) if it is  an edge of type~$1$ (type~$2$, respectively) in $G_{n}$.  $d(A)$ denotes the number of edges of  type $1$ in $A(1, 2)$. Suppose that $e= xy$ is an edge of type $1$ in $H_{n}$. Two facial $3$-cycles $xyz$ and $xyw$ contain the edge $e$. If $A(w) = A(z)$ then $e$ is called \textit{$A$-singular} (shortly \textit{singular}). If $A(w) \neq A(z)$, then $e$ is called \textit{$A$-nonsingular} (see Fisk~\cite{flobar10} and Mohar \cite{flobar14}). Let $p(A)$ denote the set of all vertices of $N_{b}$ coloured~$1$ by $A$.

\begin{lemma} \label{lemma2.1}
For each $A \in {\cal C}_{4}(H_{n})$ there is $B \in {\cal C}_{4}(H_{n })$ with $|p(B)| \leqslant 2$ such that $A$ and $B$ are equivalent up to $3\lfloor {\frac{n}{2}}\rfloor - 3p(A)$ Kempe changes and
no vertex of $p(B)$ is a single Kempe chain,
 \end{lemma} 

 \begin{PrfFact} Let $A$ be a colouring of ${\cal C}_{4}(\tilde{G}_{n })$ such that $|p(A)| > 2$. It suffices to find a colouring $B \in {\cal C}_{4}(H_{n})$ such that $|p(B)| < |p(A)|$ and $B$ is equivalent to $A$ up to $3$ Kempe changes.

  Since $p(A) > 2$, one of the following cases occurs:
   \begin{enumerate}
  \item[$(1)$] there is a vertex of $p(A)$ which is a single Kempe chain,
\item[$(2)$] there are two vertices of  $p(A)$  (say $x_1$ and $x_2$) each of which is incident with exactly one $A$-nonsingular edge,
 \item[$(3)$] there exists a path $\beta \subset N_{b}$ connecting two vertices of $p(A)$ each of which is incident with two $A$-nonsingular edges and no inner vertex of $\beta$ belongs to $p(A)$.
 \end{enumerate}

Case $(1)$ By a  \textit{trivial} Kempe change (involving only one vertex) we obtain a $4$-colouring $B$ of $H_{n }$ such that $|p(B)| < |p(A)|$.

Case $(2)$. Let $x_1$, $x_2$, $x_3 \in p(A)$ and suppose that $y_{i}z_{i}$ is an edge of $N_{a}$ with both end-vertices adjacent to $x_{i}$, for $i = 1$, $2$, $3$. Let $x_{i}y_{i}$ be the only one $A$-nonsingular edge incident with $x_i$ and suppose that $w_{i} \in N_{b}$ is adjacent both to $x_i$ and $y_i$,  for $i = 1$, $2$. We choose a Kempe chain $\xi_{1}$ ($\xi_{2}$) containing the vertex $w_{1}$  ($w_{2}$) and the vertex coloured $A(z_{1})$ ($A(z_{2})$, respectively). 

Assume first that $\xi_{1}$ does not contain $z_1$. If we switch colours of  $\xi_{1}$ we obtain a colouring $B_{1} \in {\cal C}_{4}(H_{n })$ equivalent to $A$ such that $B_{1}(w_1) = B_{1}(z_1)$. Hence,  the edge $x_{1}y_{1}$ is $B_{1}$-singular.  Since  $x_{1}z_{1}$ and $x_{1}y_{1}$ are $B_1$-singular, $\{x_{1}\}$ is a single Kempe chain. Hence, by condition $(1)$, there exists a colouring $B \in {\cal C}_{4}(H_{n })$ with $|p(B)| < |p(A)|$. 

Assume now that $\xi_{1}$ contains $z_1$. Then, $y_{3}z_{3}$ is an edge of $\xi_{1}$. Thus, $\{A(w_1), A(z_1)\} = \{A(y_3), A(z_3)\}$. We prove that  $y_{1}z_{1}$ or $y_{3}z_{3}$ is not an edge of $\xi_{2}$. Namely,  if $y_{3}z_{3} \in \xi_{2}$, then $\{A(w_2), A(z_2)\} = \{A(y_{3}), A(z_3)\}$. Hence, 
\[
\{A(w_2), A(z_2)\} = \{A(w_1), A(z_1)\} \neq \{A(y_{1}), A(z_1)\}.
\]
 Similarly, if $y_{1}z_{1} \in \xi_{2}$, then $\{A(w_2), A(z_2)\} =  \{A(y_{1}), A(z_1)\}$. Hence,
 \[
 \{A(w_2), A(z_2)\} \neq \{A(w_{1}), A(z_1)\} = \{A(y_3), A(z_3)\}.
 \]
Hence, $\xi_{2}$ does not contain $z_2$. If we switch colours on $\xi_{2}$ we obtain a colouring $B_{2} \in {\cal C}_{4}(H_{n })$ equivalent to $A$ such that $B_{2}(w_2) = B_{2}(z_2)$. Hence,  the edge $x_{2}y_{2}$ is $B_{2}$-singular. Since  $x_{2}z_{2}$ and $x_{2}y_{2}$ are $B_2$-singular, $\{x_2\}$ is a single Kempe chain. Hence, by condition $(1)$, there is a colouring $B \in {\cal C}_{4}(H_{n })$ with $|p(B)| < |p(A)|$. 

Case $(3)$. Assume that vertices $u, w \in p(A)$ are end-vertices of the path~$\beta$ satisfying condition $(3)$. Since $u$ and $w$ are both incident with two $A$-nonsingular edges, there exists a pair of $A$-nonsingular edges $ux$ and $wy$ of type $1$ such that $x, y$ are coloured the same by $A$ (say $A(x) = A(y) = i$, for some $i\in \{2, 3, 4\}$). Then, this pair splits the vertex set of $A(j, k)$ into two parts, where $\{j, k\} = \{2, 3, 4\} \setminus \{i\}$. One part of them is a Kempe chain because $\beta$ has no inner vertex belonging to $p(A)$. If we switch colours on the Kempe chain we obtain a $4$-colouring $B'$ of  $H_{n }$ equivalent to $A$ such that the edges $ux$ and $wy$ are $B'$-singular and the other edges remain singular or nonsingular. Hence, by condition $(2)$, there is a colouring $B \in {\cal C}_{4}(H_{n })$ with $|p(B)| < |p(A)|$. 
  \end{PrfFact}
  
 \begin{lemma} \label{lemma2.2}
Let $A \in {\cal C}_{4}(H_{n })$ be such  that $p(A) \leqslant 2$ and no vertex of $p(A)$ is a single Kempe chain. There is a colouring $B  \in {\cal C}_{4}(H_{n })$ such that $d(B) = p(B) \leqslant 2$ and $B$ is equivalent to $A$ up to $3p(A)$ Kempe changes. 
\end{lemma} 

\begin{PrfFact}
Let $A \in {\cal C}_{4}(H_{n})$ be such  that $p(A) \leqslant 2$ and no vertex of $p(A)$ is a single Kempe chain. 

Assume first that  $A(1, 2)$ contains only one maximal path (say $\xi$) contained in $G_{n} \setminus \{a, b\}$ (of length at least $1$). Since $p(A) \leqslant 2$, $\xi$ is of length at most $4$. It is sufficient to consider the following cases: 
\begin{enumerate}
 \item[$(a_{1})$]  $\xi$ is a path of type $2$,
  \item[$(a_{2})$] $\xi$ contains exactly two edges of type $1$,
 \item[$(a_{3})$] $\xi$ contains only one edge of type $1$.
  \end{enumerate}

\nowyiii

Case $(a_1)$. If $\xi$ is a path (of type $2$) of length $1$ or $3$, then $A(3, 4)$ is a cycle of odd order which is impossible. Hence, $\xi$ is a path of  length $2$ or $4$.

Let $\xi = x_{1}x_{2}x_{3}$ be a path of type $2$ in $A(1, 2)$ ($p(A) = 1$) and suppose that $y  \in N_{a}$ is adjacent to $x_{2}$ and $x_{3}$ (see Fig.~3).  Let $x_2$ be coloured $1$ by~$A$. If we switch colours on $\xi$ we obtain a $4$-colouring $B_{1}$ of $H_{n}$ such that~$x_2$ is coloured~$2$. Then, the edge $x_{2}y$ is a Kempe chain. If we switch colours on $x_{2}y$ we obtain a $4$-colouring $B_2$ of $H_{n}$ such that $\{x_1\}$ is a single Kempe chain. Now, we may change the colour of $x_1$ to obtain a $4$-colouring $B$ of $H_{n}$  equivalent to $A$ up to $3$ Kempe changes with $d(B) = p(B) = 1$. (The same proof is valid when $\xi$ is a path of length $4$. Then, we obtain a $4$-colouring $B$ of $H_{n}$  equivalent to $A$ with $d(B) = p(B) = 2$).%

\nowyiv 

Case $(a_2)$. Notice that $\xi$ is a path of length $4$. Let  $\xi = x_{1}\ldots x_{5}$ be a path  in $A(1, 2)$ and suppose that $x_{2}$, $x_{3}$, $x_{4}\in N_{b}$ and $x_{1}$, $x_{5}\in N_{a}$  (see Fig.~4). Then vertices $x_{2}$ and $x_{4}$ are coloured $1$ by $A$. Notice that there exists an edge $yz$ of type $2$ which is a component of $A(3, 4)$ such that $y$ is adjacent both to $x_2$ and $x_3$. If we switch colours of $yz$ we obtain  a $4$-colouring $B_1$ of $H_{n}$ such that vertices $\{x_2\}$ and $\{x_4\}$ are single Kempe chains. Hence, we may change the colours of $x_2$ and $x_4$ to obtain a $4$-colouring $B$ of  $H_{n}$ equivalent to $A$ with $d(B) = p(B) = 0$. 

\nowyv 

Case $(a_3)$. Certainly, if $\xi$ is an edge of type $1$, then lemma holds. If $\xi$ is a path of length $2$ or $4$ containing only one edge of type $1$, then $A(3, 4)$ is a path of odd order. Hence, $p(A)$ contains a vertex which is a single Kempe chain which is impossible. 

Let  $\xi = x_{1}x_{2}x_{3}x_{4}$ be a path in $A(1, 2)$ such that $x_{3}x_{4}$ is an edge of type~$1$ (see Fig.~5). Then, $x_1$ and $x_3$ are coloured $1$ by $A$. Let $y \in N_{a}$ be adjacent both to $x_1$ and $x_2$. Notice that the edge $x_{2}y$ is a Kempe chain. If we switch colours on $x_{2}y$ we obtain  a $4$-colouring $B$ of $H_{n}$ equivalent to $A$ with $d(B) = p(B) = 2$. 

Assume now that $A(1, 2)$ contains two disjoint maximal paths (say $\gamma$ and~$\delta$) contained in $G_{n} \setminus \{a, b\}$  (of lengths at least $1$). Since $p(A) \leqslant 2$, $\gamma$ and $\delta$ are both of length at most $2$. It suffices to consider the following cases:
\begin{enumerate}
 \item[$(b_{1})$] $\gamma$ and $\delta$ are paths of type $2$, 
 \item[$(b_{2})$] $\gamma$ is a path of type $2$, $\delta$ is a path of type $1$ and they are of the same length,
 \item[$(b_{3})$] $\gamma$ is a path of type $2$ and $\delta$ is an edge of type $1$, 
 \item[$(b_{4})$] $\gamma$ is a path of  type $1$ and $\delta$ is an edge of type $2$,
 \item[$(b_{5})$] $\gamma$ and $\delta$ are paths of type $1$.
 \end{enumerate}
  
\nowyvi 

Case $(b_{1})$. If  $\gamma$ is a path of  length $2$ and type $2$ and $\delta$ is and edge of type~$2$, then $A(3, 4)$ is a cycle of odd order which is impossible. Hence $\gamma$ and $\delta$ are edges or they are paths of length $2$.

Let $\gamma = x_{1}x_{2}$ and $\delta = y_{1}y_{2}$ be edges of type $2$ in $A(1, 2)$ both clockwise oriented in $N_{b}$. Suppose that $a_{1}\ldots a_{5}$ (or $b_{1}\ldots b_{5}$) is a path in $A(3, 4)$ induced by neighbours of $\{x_{1}, x_{2}\}$ ($\{y_{1}, y_{2}\}$) such that the path $a_{2}a_{3}a_{4}$ ($b_{2}b_{3}b_{4}$, respectively) is clockwise oriented in $N_{a}$. Notice that $a_3$ and $b_3$ are ends of a path of type $1$ of even order contained in $A(3, 4)$ (see Fig.~6). Then vertices $a_{1}, a_{3}, b_{4}$ have the same colour (say $A(a_1) = 3$). If $x_{1}$ and $y_{1}$ are coloured $3$ by~$A$, we switch colours on $y_{1}y_{2}$ to obtain a $4$-colouring $B_{1}$ of $H_{n}$ such that $x_{1}$ and $y_{2}$ are coloured $3$ by $B_{1}$. Hence, the pair of edges $(x_{1}a_{3}, y_{2}b_{4})$  splits the vertex set of $B_{1}(2, 4)$ into two Kempe chains. If we switch colours on vertices of the Kempe chain connecting $x_2$ and $b_3$ we obtain a $4$-colouring~$B_2$ of $H_{n}$ such that $\{x_1\}$ is a single Kempe chain. Now, we may change the colour of~$x_1$ to obtain a $4$-colouring $B$ of $H_{n}$ equivalent to $A$ with $d(B) = p(B) = 1$.  (The same proof is valid when $\gamma$ and $\delta$ are paths of type $2$ and lengths $2$. Then $p(A) = 2$. We obtain a $4$-colouring $B$ of $H_{n}$  equivalent to $A$ up to $6$ Kempe changes with $d(B) = p(B) = 2$).

Case $(b_{2})$.  If $\gamma$ and $\delta$ are paths of length $2$, $\gamma$ is of type $2$ and $\delta$ is of type~$1$, then $A(3, 4)$ is a path of odd order. Hence, $p(A)$ contains a vertex which is a single Kempe chain which is impossible.  

Let $\gamma = x_{1}x_{2}$ and $\delta = y_{1}y_{2}$ be edges of type $2$ and type $1$ in $A(1, 2)$ such that the edge $x_{1}x_{2}$ is clockwise oriented in $N_{b}$. Let $a_{1}\ldots a_{5}$ be a path in $A(3, 4)$ induced by neighbours of $\{x_{1}, x_{2}\}$ such that $a_{2}a_{3}a_{4}$ is clockwise oriented in $N_{a}$  (see Fig.~7). Suppose that $b_{1} \in N_{b}$ ($b_{2} \in N_{a}$) is adjacent both to $y_{1}$ and $y_2$ and the edge $b_{1}y_{1}$ is clockwise oriented in $N_{b}$. Notice that $a_3$ and $b_1$ ($a_3$ and~$b_2$) are ends of a path of type $1$ and odd order contained in $A(3, 4)$. Then  $a_1$, $a_3$, $b_1$ and $b_2$ have the same colour (say $A(a_1) = 3$). If $x_{2}$ is coloured $3$ by $A$, we switch colours on $x_{1}x_{2}$ to obtain a $4$-colouring $B_{1}$ of $H_{n}$ such that $x_{1}$ and $y_{1}$ are coloured $3$ by $B_{1}$.  Therefore, the pair of edges $(x_{1}a_{3}, y_{1}b_{2})$ splits the vertex set of $A(2,4)$ into two Kempe chains. If we switch colours on the Kempe chain connecting vertices $x_2$ and $y_2$ we obtain a $4$-colouring $B_2$ of $H_{n}$ such that $\{x_1\}$ and $\{y_1\}$ are single Kempe chains. Now, we may change the colour of $x_1$ and of $y_1$ to obtain a $4$-colouring $B$ of $H_{n}$ equivalent to $A$ with $d(B) = p(B) = 0$. 

\nowyvii

Case $(b_{3})$. The same proof as above is valid when $\gamma$ is a path of length $2$ and type $2$ and  $\delta$ is an edge of type $1$. Then, we obtain a $4$-colouring $B$ of $H_{n}$  equivalent to $A$ with $d(B) = p(B) = 2$.

\nowyviii

Case $(b_{4})$. Let  $\gamma = x_{1}x_{2}x_{3}$ be a path of type $1$  such that $x_{2}x_{3}$ is an edge of type $2$ clockwise oriented in $N_{b}$. Suppose that $\delta = y_{1}y_{2}$  is an edge of type $2$ in $A(1, 2)$ clockwise oriented in $N_{a}$.  Let $a_{1}\ldots a_{5}$ be a path in $A(3, 4)$ induced by neighbours of $\{y_{1}, y_{2}\}$ such that  $a_{2}a_{3}a_{4}$ is clockwise oriented in $N_{a}$ (see Fig.~8).  Suppose that $b_{1} \in N_{a}$ is adjacent both to $x_{2}$ and $x_3$.   Notice that $b_{1}$ and $a_{3}$ are ends of a path of type $1$ contained in $A(3, 4)$.  Since it is of even order,  $b_1$ and $a_4$ have the same colour (say $A(b_1) = 3$). If $x_{2}$ and $y_{1}$ are coloured $3$ by $A$, we switch colours on $y_{1}y_{2}$ to obtain a $4$-colouring $B_{1}$ of $H_{n}$ such that $x_{2}$ and $y_{2}$ are coloured $3$ by $B_{1}$.   Therefore, the pair of edges ($x_{2}b_{1}, y_{2}a_{4}$) splits the vertex set of $B_{1}(2, 4)$ into two Kempe chains. If we switch colours on  the Kempe chain containing vertices $x_3$ and $a_3$ we obtain a $4$-colouring $B$ of $H_{n}$  equivalent to $A$ with $d(B) = p(B) = 2$. 

\nowyix

Case $(b_{5})$. Certainly, if $\gamma$ and $\delta$ are edges of type $1$, then lemma holds. 

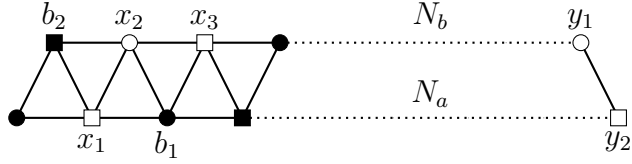
Let $\gamma = x_{1}x_{2}x_{3}$ be a path of type $1$ and $\delta$ be an edge of type $1$ in $A(1, 2)$ (see Fig.~9). Suppose that $x_{1}x_{2}$ is an edge of type $1$ and $b_{1} \in N_{a}$ ($b_2 \in N_{b}$) is adjacent both to $x_{1}$ and $x_2$. Sine $x_2$  is coloured $1$ by $A$, the pair of edges $(x_{2}x_{1}, y_{1}y_{2})$ splits the vertex set of $A(3, 4)$ into two Kempe chains. If we switch colours on one of them we obtain a $4$-colouring $B_1$ of $H_{n}$ such that $B_{1}(b_1) = B_{1}(b_2)$. Hence, $\{x_{2}\}$ is a single Kempe chain. By trivial Kempe change we obtain a $4$-colouring $B$ of $H_{n }$  equivalent to $A$ with $d(B) = p(B) = 1$. (The same proof is valid when $\gamma$ and $\delta$ are paths of type $1$ and lengths $2$. Then, we obtain a $4$-colouring $B$ of $H_{n}$  equivalent to $A$ with $d(B) = p(B) = 0$).
 \end{PrfFact}
 
\begin{lemma} \label{lemma2.3}
Let $A  \in  {\cal C}_{4}(H_{n})$ with $p(A) = d(A) = 1$. 
\begin{enumerate}
 \item[$(a)$]  If $e_{1}, e_{2}$ are edges of type $1$ with a common vertex belonging to $N_{b}$ and $e_{1}$ is the edge of $A(1, 2)$, then there exists a $4$-colouring $B$ of $H_{n}$ such that $A, B$ are equal, $p(B) = d(B) = 1$ and $e_{2}$ is the edge of $B(1, 2)$.
 \item[$(b)$] If $e, f$ are consecutive parallel edges of type $1$  and $e$ is the edge of $A(1,2)$, then there exists a $4$-colouring $B$ of $H_{n}$  such that  $A, B$ are equivalent up to $3$ Kempe changes, $p(B) = d(B) = 1$ and $f$ is the edge of $B(1, 2)$. 
\end{enumerate}
\end{lemma}

\begin{PrfFact}
 $(a)$. Let $e_{1} = a_{3}b_{2}$, $e_{2} = a_{3}b_{3}$ be edges of type $1$ and suppose that  $e_{1}$  is the edge of $A(1, 2)$, $a_{3} \in N_{b}$ and $b_{3}$ is coloured $4$ (see Fig..10) Since $p(A) = d(A) = 1$, $A(3, 4)$ is a path of type $1$ of even order. Hence, $e_1$ is nonsingular. Therefore, $A(2, 4)$ is a cycle (because $|p(A)| = 1$). If we switch colours on $A(2, 4)$ we obtain a $4$-colouring $B$ which is equal to $A$, $p(B) = d(B) = 1$ and $e_{2}$ is the edge  of $B(1, 2)$.
 
\nowyx 

Proof $(b)$.  Since $n \geqslant 5$, there exist consecutive parallel edges $a_{1}b_{1}, \ldots, a_{5}b_{5}$ of type $1$  (disjoint in pairs) such that $a_{i}\in N_{b}$, for $i = 1$, \ldots, $5$, and  $b_1$ is adjacent both to $a_1$ and $a_2$. By condition $(b)$ we may assume that $e = a_{3}b_{2} \in A(1, 2)$ and $f = a_{4}b_{3}$. Since $p(A) = d(A) = 1$, $A(3, 4)$ is a path of type $1$ of even order. Hence, $e$ is nonsingular. We may assume that $a_{2}$ and $a_{4}$ are coloured $3$ and $b_3$ is coloured $4$ by $A$ (see Fig.~10).  Since $p(A) = \{a_{3}\}$, $\delta = a_{2}a_{3}a_{4}$ is a Kempe chain of $A(1, 3)$.  If we switch colours on $\delta$ we obtain a $4$-colouring $B_{1}$ of $H_{n}$. Notice that  $\gamma = a_{3}b_{3}$ is a  Kempe chain of $B_{1}(3, 4)$. If we switch colours on $\gamma$  we obtain a $4$-colouring $B_{2}$ of $H_{n}$. Notice that $a_2$ is a single Kempe chain of $B_{2}(1, 3)$. We may change colour of $a_2$ to  obtain a $4$-colouring $B_3$ of $H_{n}$. Since $p(B_{3}) = \{a_{4} \}$, $a_{3}a_{4}a_{5}$ is a Kempe path of $B_{3}(1, 4)$.  Hence, $B_{3}(2, 3)$ is a cycle (see Fig.~10). Next we switch colours on $B_{3}(2, 3)$ to obtain a $4$-colouring $B$ which is equal to $B_{3}$ such that $p(B) = d(B) = 1$ and $f = a_{4}b_{3}$ is the edge of $B(1, 2)$. Certainly, colourings $A$ and $B$ are equivalent up to $3$ Kempe changes. 
\end{PrfFact}

\begin{corollary}\label{corollary2.1}
Let $A  \in  {\cal C}_{4}(H_{n})$ such that $p(A) = d(A) = 2$ and two edges of $A(1, 2)$ are nonsingular. 
 \begin{enumerate}
\item[$(c)$] If $e_{1}, e_{2}$ are edges of type $1$ with a common vertex belonging to~$N_{b}$ and  $e_{1}$ is a Kempe chain of $A(3, 4)$, then there exists a $4$-colouring $B$ of $H_{n}$ equal to $A$ such that $p(B) = d(B) = 2$, two edges of $B(1, 2)$ are $B$-nonsingular and $e_{2}$ is a Kempe chain of $B(3, 4)$.
 \item[$(d)$] If $e, f$ are consecutive parallel edges of type $1$ and  $e$ is a Kempe chain of $A(3, 4)$, then there exists a $4$-colouring $B$ of $H_{n}$ such that $A$, $B$ are equivalent up to $3$ Kempe changes, $p(B) = d(B) = 2$, two edges of $B(1, 2)$ are $B$-nonsingular and $f$ is a Kempe chain of $B(3, 4)$.
\end{enumerate}
\end{corollary}

 \begin{PrfFact} (c)--(d). The same proof as for condition $(a)$  (or $(b)$) remains valid for $(c)$ ($(d)$, respectively).
\end{PrfFact} 

 \begin{definition}\label{definition2.1}
 If $A  \in {\cal C}_{4}(H_{n })$ and $A(1, 2)$ has no edge of type $2$, then there exists a colouring $A^{+}  \in {\cal C}_{4}(H_{n})$  defined in the following way: if $v \in N_{b}$, $A(v) = 1$ and $v$ is not a vertex of $A (1, 2)$, then $A^{+}(v) = 2$ and $A^{+}(w) = A(w)$ for other vertices of $H_{n}$. Notice that $d(A^{+}) = d(A)$. 
 
If $A \in {\cal C}_{4}(H_{n})$ and $A(1, 2)$ has no edge of type $2$, then there exists a colouring $A^{-}  \in {\cal C}_{4}(H_{n})$  defined in the following way: if $v \in N_{b}$, $A(v) =~2$, then $A^{-}(v) = 1$ and $A^{-}(w) = A(w)$ for other vertices of $H_{n}$. Notice that $d(A^{-}) = d(A)$. 

If $A \in {\cal C}_{4}(H_{n})$ and $A(1, 2)$ has no edge of type $2$, then there exists a colouring $\overline{A} \in {\cal C}_{4}(G_{n})$  defined in the following way: $\overline{A}(b) = 2$ and $\overline{A}(w) = A^{-}(w)$ for other vertices of $H_{n}$. Notice that $d(\overline{A}) = d(A)$.

If $A \in {\cal C}_{4}(G_{n})$, then $A|(H_{n})$ denote the restriction of $A$ to $H_{n}$. 
\end{definition}

\begin{theorem} \label{theorem2.1} For every graph $H_{n}$, $n \geqslant 5$, $4$-colourings of  $H_{n}$ are all  equivalent up to
\begin{align*}
&  6\left\lfloor \frac{n}{2} \right\rfloor \ \hbox{Kempe changes}, \ \hbox{for} \  n \equiv 0\, 
\pmod3,
\\
& 9\left\lfloor \frac{n}{2} \right\rfloor \ \hbox{Kempe changes},  \ \hbox{for} \  n \equiv 2\,
\pmod3,
\\
&9\left\lfloor {\frac{n}{2}}\right\rfloor + 6\left\lfloor {\frac{n}{3}}\right\rfloor - 2\ \hbox{Kempe changes},  \ \hbox{for} \  n \equiv 1\, 
\pmod3.
\end{align*}
\end{theorem}

\begin{PrfFact}
Let $A_{0}$, $C_{0} \in {\cal C}_{4}(H_{n})$. In view of Lemma \ref{lemma2.1} and Lemma \ref{lemma2.2},  there exists colouring $A \in {\cal C}_{4}(H_{n})$ (or $C \in {\cal C}_{4}(H_{n})$) with $p(A) = d(A) \leqslant 2$ ($p(C) = d(C) \leqslant 2$) which is equivalent to $A_{0}$ ($C_{0}$, respectively) up to $3\lfloor \frac{n}{2} \rfloor$ Kempe changes.  Since $A(1, 2)$ has no edge of type $2$, there exists  the colouring $\overline{A} \in {\cal C}_{4}({G}_{n})$ with $d(\overline{A}) = d(A) \leqslant 2$. By condition $(4)$ of Lemma \ref{lemma1.3}, we obtain
\[
n \equiv 2k\, 
\pmod3
\hbox{ if and only if } d(\overline{A}) = k, \hbox{ for  $k = 0$, $1$, $2$}.
\]

If $n \equiv 0\, 
\pmod3$, then $d(\overline{A}) = 0$. Hence, $p(A) = d(A) = 0$. Then, $A = Q|H_{n}$. Similarly,  $C = Q|H_{n}$. Hence,  $A_{0}$ and $C_{0}$ are equivalent up to $6\lfloor {\frac{n}{2}}\rfloor$ Kempe changes.

If $n \equiv 2\, 
\pmod3$, then  $d(\overline{A}) = 1$. Hence, $p(A) = d(A) = 1$. Similarly, $p(C) = d(C) = 1$. By Lemma \ref{lemma2.3},  $A$ and $C$ are equivalent up to $3\lfloor {\frac{n}{2}}\rfloor$ Kempe changes. Thus, $A_{0}$ and $C_{0}$ are equivalent up to $9\lfloor {\frac{n}{2}}\rfloor$ Kempe changes. 

If $n \equiv 1\, 
\pmod$, then, $d(\overline{A}) = 2$. Hence, $p(A) = d(A) = 2$. Thus, colourings $A^{-}$ and  $A$ are equivalent up to $a(\overline{A}) - 3$ single Kempe changes (where $a(\overline{A})$ is the number of vertices coloured $1$ by $\overline{A}$). 
 If  edges of $A^{-}(1,2)$ are nonsingular, then, by  $d(A^{-}) = 2$, we may switch colours on vertices in one of two Kempe chains contained  in $A^{-}(3, 4)$ to obtain a $4$-colouring  $A^{(-, s)}$ of $H_{n}$which  edges are singular. Then, by Corollary~\ref{corollary1.2},  $\overline{A^{(-, s)}}$ and $Q_{2,e}$ are  equivalent colourings of 
 $G_{n}$ up to  $\frac {c(\overline{A})}{2}$ Kempe changes each of which switches colours on vertices in some edge of type $2$ in $G_{n}$. Hence, $\overline{A^{(-, s)}}|H_{n}$ and $Q_{2,e}| H_{n}$ are equivalent colourings of $H_{n}$ up to $\frac {c(\overline{A})}{2}$ Kempe changes.

Further, $Q_{2,e}|H_{n}$ and $(Q_{2,e}|H_{n})^{+}$ are equivalent up to $a(Q_{2, e})- 3$ single Kempe changes. If edges of $(Q_{2,e}|H_{n})^{+}(1, 2)$ are singular, then, by $d((Q_{2,e}|H_{n})^{+}) = 2$, we may switch colours on vertices in one of two Kempe chains in $(Q_{2,e}|H_{n})^{+}(3, 4)$ to obtain a $4$-colouring  $B = (Q_{2,e}|H_{n})^{(+, ns)}$ of $H_{n}$ such that edges   of $B(1, 2)$ are nonsingular.  
 
Notice that, by conditions $(3)$ and $(4)$ of  Lemma \ref{lemma1.3},  
\[
a(\overline{A}) - 3 =  \frac {c(\overline{A})}{2} = \frac{n - 4}{3}\leqslant  \left\lfloor {\frac{n}{3}}\right\rfloor - 1 \ \hbox{ and } 
\  a(Q_{2, e}) - 3 = \frac{n - 4}{3}.
\]
Therefore, 
 \[
 A \sim A^{-} \sim A^{(-, s)} =  \overline{A^{(-, s)})}|H_{n} \sim Q_{2,e}|H_{n}  \sim (Q_{2,e}|H_{n})^{+} \sim (Q_{2,e}|H_{n})^{(+, ns)}
 \]
 up to $3(\lfloor {\frac{n}{3}}\rfloor  - 1) + 2$ Kempe changes. Similarly, $C$ and $(Q_{2, f}|H_{n})^{(+, ns)}$ are equivalent up to $3 \lfloor {\frac{n}{3}}\rfloor - 1$ Kempe changes.

By Corollary \ref{corollary2.1}, $(Q_{2, e}|H_{n})^{(+, ns)}$ and $(Q_{2, f}|H_{n})^{(+, ns)}$ are equivalent up to $3\lfloor {\frac{n}{2}}\rfloor$ Kempe changes, for every edges $e$ and $f$ of type~$1$. Hence, $A_{0}$ and $C_{0}$ are equivalent up to  $9\lfloor {\frac{n}{2}}\rfloor + 6\lfloor {\frac{n}{3}}\rfloor - 2$ Kempe changes.
\end{PrfFact}

\end{document}